\documentclass[twocolumn,journal]{IEEEtran}
\usepackage[latin9]{inputenc}
\usepackage{mathtools}
\usepackage{amsmath}
\usepackage{amsthm}
\usepackage{amssymb}
\usepackage{cancel}
\usepackage{mathdots}
\usepackage{stmaryrd}
\usepackage{stackrel}
\usepackage{undertilde}
\PassOptionsToPackage{version=3}{mhchem}
\usepackage{mhchem}
\usepackage{esint}
\PassOptionsToPackage{normalem}{ulem}
\usepackage{ulem}
\usepackage[unicode=true,
 bookmarks=true,bookmarksnumbered=true,bookmarksopen=true,bookmarksopenlevel=1,
 breaklinks=true,pdfborder={0 0 0},pdfborderstyle={},backref=false,colorlinks=false]
 {hyperref}
\hypersetup{pdftitle={On Convergence of Heuristics Based on Douglas-Rachford Splitting and ADMM to Minimize Convex Functions over Nonconvex Sets},
 pdfauthor={Shuvomoy Das Gupta},
 pdfpagelayout=OneColumn, pdfnewwindow=true, pdfstartview=XYZ, plainpages=false}

\makeatletter
\theoremstyle{definition}
\newtheorem{example}{\protect\examplename}
\theoremstyle{plain}
\newtheorem{thm}{\protect\theoremname}
\theoremstyle{plain}
\newtheorem{lem}{\protect\lemmaname}
\theoremstyle{remark}
\newtheorem{rem}{\protect\remarkname}
\theoremstyle{plain}
\newtheorem{cor}{\protect\corollaryname}
\theoremstyle{plain}
\newtheorem{prop}{\protect\propositionname}

\IEEEoverridecommandlockouts
\usepackage[caption=false,font=footnotesize]{subfig}
\usepackage{tikz} 
\allowdisplaybreaks
\begingroup
\makeatletter
   \@for\theoremstyle:=definition,remark,plain\do{%
     \expandafter\g@addto@macro\csname th@\theoremstyle\endcsname{%
        \addtolength\thm@preskip\parskip
     }%
   }
\endgroup
\usepackage{parskip}

\makeatletter
\def\@IEEEsectpunct{.\ \,}
\def\paragraph{\@startsection{paragraph}{4}{\z@}{1ex plus 1ex minus 0.5ex}%
{0ex}{\normalfont\normalsize\bfseries}}
\makeatother

\makeatother

\providecommand{\corollaryname}{Corollary}
\providecommand{\examplename}{Example}
\providecommand{\lemmaname}{Lemma}
\providecommand{\propositionname}{Proposition}
\providecommand{\remarkname}{Remark}
\providecommand{\theoremname}{Theorem}

\begin{document}
\title{On Convergence of Heuristics Based on Douglas-Rachford Splitting and
ADMM to Minimize Convex Functions over Nonconvex Sets}
\author{Shuvomoy Das Gupta \thanks{S. Das Gupta is with the Research \& Technology Department, Thales
Canada, 105 Moatfield Drive, Toronto, ON, Canada. Email: \texttt{\protect\href{mailto:shuvomoy.dasgupta@thalesgroup.com}{shuvomoy.dasgupta@thalesgroup.com}}}}
\maketitle
\begin{abstract}
Recently, heuristics based on the Douglas-Rachford splitting algorithm
and the alternating direction method of multipliers (ADMM) have found
empirical success in minimizing convex functions over nonconvex sets,
but not much has been done to improve the theoretical understanding
of them. In this paper, we investigate convergence of these heuristics.
First, we characterize optimal solutions of minimization problems
involving convex cost functions over nonconvex constraint sets. We
show that these optimal solutions are related to the fixed point set
of the underlying nonconvex Douglas-Rachford operator. Next, we establish
sufficient conditions under which the Douglas-Rachford splitting heuristic
either converges to a point or its cluster points form a nonempty
compact connected set. In the case where the heuristic converges to
a point, we establish sufficient conditions for that point to be an
optimal solution. Then, we discuss how the ADMM heuristic can be constructed
from the Douglas-Rachford splitting algorithm. We show that, unlike
in the convex case, the algorithms in our nonconvex setup are not
equivalent to each other and have a rather involved relationship between
them. Finally, we comment on convergence of the ADMM heuristic and
compare it with the Douglas-Rachford splitting heuristic. 
\end{abstract}


\section{Introduction}

In this paper, we study convergence of heuristics based on the Douglas-Rachford
splitting algorithm and the alternating direction method of multipliers
(ADMM) for minimization of convex functions over nonconvex sets. Such
optimization problems can be described as 
\begin{equation}
\begin{array}{ll}
\textup{minimize} & f\left(x\right)\\
\textup{subject to} & x\in\mathcal{C},
\end{array}\tag{OPT}\label{eq:original_problem}
\end{equation}
where $x\in\mathbf{R}^{n}$ is the decision variable. The constraint
set $\mathcal{C}\subseteq\mathbf{R}^{n}$ is nonempty and compact
(closed and bounded), but it is not necessarily convex. The cost function
$f:\mathbf{R}^{n}\to\mathbf{R}\cup\left\{ +\infty\right\} $ is CPC.
This means that $f$ is (i) \uline{c}onvex, (ii) \uline{p}roper,
\emph{i.e.,} its domain $\mathop{{\bf dom}}f=\{x\in\mathbf{R}^{n}\mid f(x)<\infty\}$
is nonempty, and (iii) \uline{c}losed (or lower-semicontinuous),
\emph{i.e.,} its epigraph $\mathop{{\bf epi}}f=\left\{ (x,\xi)\in\mathbf{R}^{n}\times\mathbf{R}\mid f(x)\leq\xi\right\} $
is a closed set. The constraint set $\mathcal{C}$ is assumed to be
closed so that projection onto $\mathcal{C}$ is well-defined, and
it is assumed to be bounded to avoid the possibility of an unbounded
optimal solution.

We consider the following \textbf{heuristic based on the Douglas-Rachford
splitting algorithm }\cite[\S 27.2]{Bauschke2017} to solve \eqref{eq:original_problem}:
\begin{equation}
\begin{aligned}x_{n+1} & =\mathbf{prox}_{\gamma f}\left(z_{n}\right)\\
y_{n+1} & =\tilde{\mathbf{\Pi}}_{\mathcal{C}}\left(2x_{n+1}-z_{n}\right)\\
z_{n+1} & =z_{n}+y_{n+1}-x_{n+1},
\end{aligned}
\tag{NC-DRS}\label{eq:DR_splitting_orig}
\end{equation}
where $n\in\mathbf{N}$ is an iteration counter, $\tilde{\mathbf{\Pi}}_{\mathcal{C}}(x)$
is a \textbf{Euclidean projection} of $x$ onto $\mathcal{C}$ (as
$\mathcal{C}$ is not necessarily convex, there can be multiple projections
onto it from a point outside $\mathcal{C}$), and 
\[
\mathbf{prox}_{\gamma f}\left(x\right)=\textrm{argmin}_{y\in\mathop{{\bf dom}}f}\big(f(y)+\frac{1}{2\gamma}\|y-x\|^{2}\big)
\]
is the \textbf{proximal operator} of $f$ at $x$ with parameter $\gamma>0$.
We also investigate the following \textbf{heuristic based on ADMM}
(also known as NC-ADMM \cite[\S 3.1]{Diamond2018}) to solve \eqref{eq:original_problem}:
\begin{equation}
\begin{aligned}x_{n+1} & =\mathbf{prox}_{\gamma f}\left(y_{n}-z_{n}\right)\\
y_{n+1} & =\tilde{\mathbf{\Pi}}_{\mathcal{C}}\left(x_{n+1}+z_{n}\right)\\
z_{n+1} & =z_{n}-y_{n+1}+x_{n+1},
\end{aligned}
\tag{NC-ADMM}\label{eq:ADMM_orig-1}
\end{equation}
where $n\in\mathbf{N}$ is an iteration counter, and $\gamma>0$.
Note that both heuristics consist of the same subroutines, but different
inputs are fed into them. So, any software package that implements
one of the heuristics can be easily modified to implement the other. 

If the set $\mathcal{C}$ is convex, then the iterates $x_{n},y_{n}$
in both \eqref{eq:DR_splitting_orig} and \eqref{eq:ADMM_orig-1}
converge to an optimal solution for any initial point \cite[Corollary 27.4]{Bauschke2017},
\cite{Boyd2011}. The convergence conditions for the nonconvex case,
studied in this paper, are far more complicated.

\paragraph*{Motivation}

This paper is motivated by the recent success of ADMM in solving nonconvex
problems. ADMM, which is a special case of the Douglas-Rachford splitting
algorithm in a convex setup, was originally designed to solve convex
optimization problems \cite{Boyd2011}. However, since the idea of
implementing this algorithm as a general purpose heuristics to solve
nonconvex optimization problems was introduced in \cite{Boyd2011},
\eqref{eq:ADMM_orig-1} has been applied successfully to minimization
of convex functions over nonconvex sets \cite{Diamond2018,Takapoui2017a,Takapoui2017,takapoui2016simple},
and it has been implemented recently in the Python package NCVX\textemdash an
extension of CVXPY\textemdash to formulate and solve problems of the
form \eqref{eq:original_problem} \cite{diamond2016cvxpy}. In these
works, the nonconvex projection step of \eqref{eq:ADMM_orig-1}, if
computationally too costly, is replaced with a tractable ``approximate''
projection onto the nonconvex set, \emph{e.g.,} rounding for Boolean
variables; yet it finds approximate solutions to a wide variety of
nonconvex problems effectively. In spite of the empirical success,
not much has been done to improve the theoretical understanding of
such heuristics. Some recent progress has been made on understanding
convergence of ADMM for specialized nonconvex setups, such as (i)
minimizing a nonconvex function over an affine set \cite{wang2015global},
and (ii) minimizing the sum of a smooth function with a bounded Hessian
and a nonsmooth function with an easy to compute proximal mapping
\cite{li2015global}. However, these works are not applicable to
\eqref{eq:DR_splitting_orig} and \eqref{eq:ADMM_orig-1}, which has
motivated us to investigate the convergence properties of these heuristics. 

\paragraph*{Contributions}

Our contributions in this paper are as follows. First, we characterize
global minimizers of \eqref{eq:original_problem} and show that they
are related to the fixed point set of the underlying Douglas-Rachford
operator constructed from \eqref{eq:DR_splitting_orig}. Then, we
establish conditions under which \eqref{eq:DR_splitting_orig} either
converges to a point (not necessarily an optimal solution) or its
cluster points form a nonempty compact connected set. In the case
where the heuristic converges to a point, we provide sufficient conditions
for that point to be an optimal solution. Then, we investigate the
relationship between \eqref{eq:DR_splitting_orig} and \eqref{eq:ADMM_orig-1}.
For a convex optimization problem, ADMM is the Douglas-Rachford algorithm
splitting applied to the dual problem \cite{Gabay1983}, but their
relationship is more involved in our nonconvex setup. Applying the
Douglas-Rachford splitting algorithm to the convex dual of \eqref{eq:original_problem}
results in a relaxed version of \eqref{eq:ADMM_orig-1}, where the
projection is onto the convex hull of $\mathcal{C}$. We show that
this relaxed algorithm finds a minimizer of $f$ over the convex hull
of $\mathcal{C}$, and by restricting its projection step onto the
original constraint set $\mathcal{C}$, we arrive at \eqref{eq:ADMM_orig-1}.
The construction procedure also explains why, when compared with exact
solvers, \eqref{eq:ADMM_orig-1} often achieves lower objective values
in many numerical experiments performed in \cite{Diamond2018,Takapoui2017a,Takapoui2017,takapoui2016simple}.
We comment on the convergence properties of \eqref{eq:ADMM_orig-1}
and compare it with \eqref{eq:DR_splitting_orig}. To the best of
our knowledge, we are not aware of similar results in the existing
literature.

\paragraph*{Notation and notions}

We denote the sets of real numbers and natural numbers by $\mathbf{R}$
and $\mathbf{N}$, respectively. Furthermore, $\overline{\mathbf{R}}=\mathbf{R}\cup\{\infty\}$
denotes the extended real line. The set of real column vectors of
length $n$ is denoted by $\mathbf{R}^{n}$. Depending on the context,
$0$ may be a scalar or a column vector of zeros. The $n\times n$
identity matrix is denoted by $I_{n}$. The standard Euclidean norm
is denoted by $\|\cdot\|$. We use $\left\langle \cdot\mid\cdot\right\rangle $
as the inner product in the Euclidean space. Let $\mathcal{X},\mathcal{Y}$
be two nonempty subsets of $\mathbf{R}^{n}$, and let $z\in\mathbf{R}^{n}$.
Then, $\mathcal{X}+\mathcal{Y}=\{x+y\mid x\in\mathcal{X},y\in\mathcal{Y}\}$,
$\mathcal{X}-\mathcal{Y}=\{x-y\mid x\in\mathcal{X},y\in\mathcal{Y}\}$,
$z+\mathcal{X}=\{z\}+\mathcal{X}$, and $\mathcal{X}-z=\mathcal{X}-\{z\}$.
If one of the sets is empty, then the resultant addition or subtraction
is an empty set, \emph{i.e.,} $\mathcal{X}+\emptyset=\emptyset$.
Finally, the\textbf{ indicator function} of a nonempty set $\mathcal{X}\subseteq\mathbf{R}^{n}$,
denoted by $\delta_{\mathcal{X}}$, is defined as
\[
\delta_{\mathcal{X}}(x)=\begin{cases}
0, & \textrm{if}\;x\in\mathcal{X}\\
\infty, & \textrm{if}\;x\notin\mathcal{X}.
\end{cases}
\]

Using indicator function, \eqref{eq:original_problem} can be expressed
as $\textrm{minimize}\;f(x)+\delta_{\mathcal{C}}(x)$.

\section{Background on monotone operator theory\label{sec:Background-and-preliminaries}}

In this section, we present some definitions and preliminary results
on monotone operator theory and relate them to our setup. First, in
\S\ref{subsec:Operator-theoretic-notion}, we briefly review the
essential operator theoretic notions and provide examples that relate
these concepts to \eqref{eq:original_problem}. In \S\ref{subsec:Nonexpansive,-firmly-nonexpansiv},
we review nonexpansiveness and its several variants for an operator.
These concepts are essential for characterizing (i) the operators
$\tilde{\mathbf{\Pi}}_{\mathcal{C}}$ and $\mathbf{prox}_{\gamma f}$
and (ii) the fixed point sets of the underlying operators of \eqref{eq:DR_splitting_orig}
and \eqref{eq:ADMM_orig-1}. Finally, in \S\ref{subsec:Resolvent-and-reflected},
we introduce resolvent and reflected resolvent of an operator to provide
additional characterizations of $\tilde{\mathbf{\Pi}}_{\mathcal{C}}$
and $\mathbf{prox}_{\gamma f}$. 

\subsection{Operator theoretic notions \label{subsec:Operator-theoretic-notion}}

A\textbf{ set-valued operator}\emph{ }$T:\mathbf{R}^{n}\rightrightarrows\mathbf{R}^{n}$
maps each element in $\mathbf{R}^{n}$ to a set in $\mathbf{R}^{n}$;
its\emph{ }\textbf{domain} is defined as $\mathop{{\bf dom}}T=\left\{ x\in\mathbf{R}^{n}\mid T(x)\neq\emptyset\right\} $,
its \textbf{range} is defined as $\mathop{{\bf ran}}T=\bigcup_{x\in\mathbf{R}^{n}}T(x)$,
and it is completely completely characterized by its \textbf{graph}:
$\mathop{{\bf gra}}T=\{(x,u)\mid u\in Tx\}.$ Furthermore, the \textbf{zero
set} of $T$ is defined as $\mathop{{\bf zer}}T=\left\{ x\mid0\in A\left(x\right)\right\} $,
and the\textbf{ fixed point set} of $T$ is defined as $\mathop{{\bf fix}}T=\left\{ x\mid Tx\ni x\right\} $.
The operator $2T-I_{n}$ is called the \textbf{reflection of} $T$.
\textbf{Inverse} of $T$, denoted by $T^{-1}$, is defined through
its graph: $\mathop{{\bf gra}}T^{-1}=\left\{ (u,x)\mid(x,u)\in\mathop{{\bf gra}}T\right\} $,
so $x\in T(u)\Leftrightarrow u\in T^{-1}(x)$. For every $x,$ \textbf{addition
of two operators} $T_{1},T_{2}:\mathbf{R}^{n}\rightrightarrows\mathbf{R}^{n}$,
denoted by $T_{1}+T_{2}$, is defined as $\left(T_{1}+T_{2}\right)\left(x\right)=T_{1}\left(x\right)+T_{2}\left(x\right)$
(subtraction is defined analogously), and \textbf{composition} of
these operators, denoted by $T_{1}T_{2},$ is defined as $T_{1}T_{2}\left(x\right)=T_{1}\left(T_{2}\left(x\right)\right)$;
note that order matters for composition. Also, if $\mathcal{X}\subseteq\mathbf{R}^{n}$
is a nonempty set, then $T(\mathcal{X})=\bigcup_{x\in\mathcal{X}}T(x)$. 

A set-valued operator $T:\mathbf{R}^{n}\rightrightarrows\mathbf{R}^{n}$
is \textbf{monotone} if, for every $\left(x,u\right),(y,v)\in\mathop{{\bf gra}}T$,
it satisfies $\left\langle x-y\mid u-v\right\rangle \geq0$. A monotone
operator $T$ is \textbf{maximally monotone} if $\mathop{{\bf gra}}T$
is not properly contained by the graph of any other monotone operator. 

Finally, a\textbf{ single-valued operator} $T:\mathcal{D}\to\mathbf{R}^{n}$
is a special type of set-valued operator, which maps every $x$ in
its domain $\mathcal{D}\subseteq\mathbf{R}^{n}$ to a singleton $T\left(x\right)$
in $\mathbf{R}^{n}$.
\begin{example}[\textbf{projection operator}]
 \label{exa:(Projection-operator)}Recall that $\tilde{\mathbf{\Pi}}_{\mathcal{C}}(x)$
is a\textbf{ }Euclidean projection of $x$ onto the constraint set
$\mathcal{C}$. The operator $\tilde{\mathbf{\Pi}}_{\mathcal{C}}$
is single-valued. The set of all projections onto $\mathcal{C}$,
denoted by $\mathbf{\Pi}_{\mathcal{C}}$, is the \textbf{set-valued
projection operator} onto $\mathcal{C}$, and it is defined as $\mathbf{\Pi}_{\mathcal{C}}\left(x\right)=\textrm{argmin}_{y\in\mathcal{C}}\|x-y\|^{2}$.
Clearly, $\tilde{\mathbf{\Pi}}_{\mathcal{C}}\left(x\right)\subseteq\mathbf{\Pi}_{\mathcal{C}}\left(x\right)$
for every $x$. Both $\mathbf{\Pi}_{\mathcal{C}}$ and $\tilde{\mathbf{\Pi}}_{\mathcal{C}}$
are monotone operators, but not necessarily maximally monotone \cite[Example 20.12]{Bauschke2017}.
The projection operator onto a nonempty closed convex set, however,
is maximally monotone \cite[Example 20.12, Corollary 20.27, and Proposition 4.8]{Bauschke2017}.
\end{example}
\begin{example}[\textbf{subdifferential operator}]
\textbf{ }\label{exa:Subdifferential-operator}For every proper function
$g:\mathbf{R}^{n}\to\overline{\mathbf{R}}$, its subdifferential operator
is the set-valued operator $\partial g:\mathbf{R}^{n}\rightrightarrows\mathbf{R}^{n}$,
which is defined as
\begin{equation}
\partial g\left(x\right)=\left\{ u\in\mathbf{R}^{n}\mid\left(\forall y\in\mathbf{R}^{n}\right)\,g\left(y\right)\geq g\left(x\right)+\left\langle u\mid y-x\right\rangle \right\} .\label{eq:subdifferential}
\end{equation}

A vector $u\in\partial g\left(x\right)$ is called a \textbf{subgradient}
of $g$ at $x$. The subdifferential operator of a proper function
is monotone, hence $\partial\delta_{\mathcal{C}}$ is monotone \cite[Example 20.3]{Bauschke2017}.
On the other hand, the subdifferential operator of a CPC function
is maximally monotone, thus $\partial f$ is maximally monotone \cite[Theorem 20.40]{Bauschke2017}.
The following result regarding the subdifferential operator plays
a key role in characterizing global minimizers of \eqref{eq:original_problem}
in $\S$\ref{sec:Convergence-analysis-for-NCDR}.
\end{example}
\begin{thm}[{\textbf{Fermat's rule }\cite[page 223]{Bauschke2017}\emph{,} \cite[\S 2.3]{subgradBoyd}}]
\textup{}\textbf{\textup{\label{thm:(Fermat's-rule)}}} The set of
all global minimizers of a proper function $g:\mathbf{R}^{n}\to\overline{\mathbf{R}}$,
denoted by $\mathop{{\rm argmin}}g$, is equal to the zero set of
its subdifferential operator $\partial g$, i.e., $\mathop{\textup{argmin}}g=\mathop{{\bf zer}}\partial g=\left\{ x\in\mathbf{R}^{n}\mid0\in\partial g(x)\right\} .$
\end{thm}
\begin{IEEEproof}
Take $x\in\mathop{{\rm argmin}}g$ which is equivalent to the statement
$\left(\forall y\in\mathbf{R}^{n}\right)\;g(y)\geq g(x)+\left\langle 0\mid y-x\right\rangle \Leftrightarrow\partial g(x)\ni0\Leftrightarrow x\in\mathop{{\bf zer}}\partial g.$ 
\end{IEEEproof}
While this simple characterization of optimality via the subdifferential
holds for every nonconvex functions, it may not be particularly useful
in practice if we cannot compute the subdifferential in an algorithmic
manner \cite[page 4]{subgradBoyd}. 

We now present a lemma regarding the subdifferential operator of the
sum of two proper functions, which is used later in $\S$\ref{sec:Convergence-analysis-for-NCDR}.
Recall that $\partial\left(g+h\right)\left(x\right)=\partial\left(g(x)+h(x)\right)$
according to our notation.
\begin{lem}[\textbf{subdifferential of sum of proper functions}]
\label{lem:subdifferential-monotone} Let $g:\mathbf{R}^{n}\to\overline{\mathbf{R}}$
and $h:\mathbf{R}^{n}\to\overline{\mathbf{R}}$ be proper functions
such that $\mathop{{\bf dom}}g\cap\mathop{{\bf dom}}h\neq\emptyset$.
Then,

\emph{(i)} the function $g+h$ is proper,

\emph{(ii)} for every $x$ in $\mathbf{R}^{n}$, we have $\partial g\left(x\right)+\partial h\left(x\right)\subseteq\partial(g+h)(x)$,
and

\emph{(iii)} both $\partial(g+h)$ and $\partial g+\partial h$ are
monotone operators.
\end{lem}
\begin{IEEEproof}
(i): By definition, $\mathop{{\bf dom}}\left(g+h\right)=\left\{ x\mid g(x)+h(x)<\infty\right\} =\mathop{{\bf dom}}g\cap\mathop{{\bf dom}}h\neq\emptyset.$
Thus, $g+h$ is proper.

(ii): Take $x\in\mathbf{R}^{n}$, and denote $u\in\partial g\left(x\right)$
and $v\in\partial h\left(x\right)$. We want to prove that $u+v\in\partial\left(g+h\right)\left(x\right)=\partial\left(g(x)+h(x)\right)$.
Using \eqref{eq:subdifferential}, we have $g\left(y\right)\geq g\left(x\right)+\left\langle u\mid y-x\right\rangle $
and $h\left(y\right)\geq h\left(x\right)+\left\langle v\mid y-x\right\rangle $
for every $y\in\mathbf{R}^{n}$. Adding the last two inequalities
we get $\left(g\left(y\right)+h\left(y\right)\right)\geq\left(g\left(x\right)+h\left(x\right)\right)+\left\langle u+v\mid y-x\right\rangle $
for every $y\in\mathbf{R}^{n}$, \emph{i.e.,} $u+v\in\partial\left(g\left(x\right)+h\left(x\right)\right)$. 

(iii): Denote $\phi:=g+h$, which is proper due to (i). Now take $(x,u),(y,v)$
in $\mathop{{\bf gra}}\partial\phi$, so we have $\phi\left(y\right)\geq\phi\left(x\right)+\left\langle u\mid y-x\right\rangle $
and $\phi(x)\geq\phi(y)+\left\langle v\mid x-y\right\rangle $ using
\eqref{eq:subdifferential}; adding these inequalities we have $0\geq\left\langle u\mid y-x\right\rangle +\left\langle v\mid x-y\right\rangle $
\emph{i.e.,} $\left\langle u-v\mid x-y\right\rangle \geq0$, so $\partial\phi=\partial\left(g+h\right)$
is a monotone operator by definition. Furthermore, both $\partial g$
and $\partial h$ are monotone, as the subdifferential operator of
a proper function is monotone \cite[Example 20.3]{Bauschke2017}.
Using also the fact that sum of two monotone operators is a monotone
operator \cite[page 351]{Bauschke2017}, we conclude that $\partial g+\partial h$
is monotone. 
\end{IEEEproof}

\subsection{Nonexpansive and firmly nonexpansive operator\label{subsec:Nonexpansive,-firmly-nonexpansiv}}

Let $T:\mathcal{D}\to\mathbf{R}^{n}$ be a single-valued operator,
where $\mathcal{D}\subseteq\mathbf{R}^{n}$ is nonempty. Then, $T$
is
\begin{enumerate}
\item \textbf{nonexpansive}\emph{ }on $\mathcal{D}$ if for every $x,y\in\mathcal{D}$
it satisfies $\left\Vert T\left(x\right)-T\left(y\right)\right\Vert \leq\|x-y\|$,
and
\item \textbf{firmly nonexpansive} on $\mathcal{D}$ if for every $x,y\in\mathcal{D}$
it satisfies $\|T\left(x\right)-T\left(y\right)\|^{2}+\left\Vert (I_{n}-T)(x)-(I_{n}-T)(y)\right\Vert ^{2}\leq\|x-y\|^{2}.$
\end{enumerate}
An operator $T:\mathcal{D}\to\mathbf{R}^{n}$ is firmly nonexpansive
on $\mathcal{D}$ if and only if its reflection operator $2T-I_{n}$
is nonexpansive \cite[Proposition 4.2]{Bauschke2017}. Furthermore,
a firmly nonexpansive operator is also nonexpansive \cite[page 59]{Bauschke2017}. 
\begin{example}[\textbf{proximal operator}]
\textbf{}\label{(proximal-operator)-The}The proximal operator of
a CPC function is both firmly nonexpansive and nonexpansive \cite[Proposition 12.27, Example 23.3]{Bauschke2017},
hence $\mathbf{prox}_{\gamma f}$ in \eqref{eq:ADMM_orig-1} and \eqref{eq:DR_splitting_orig}
is both firmly nonexpansive and nonexpansive. Furthermore, its reflection
$2\,\mathbf{prox}_{\gamma f}-I_{n}$ is nonexpansive \cite[Proposition 4.2]{Bauschke2017}. 
\end{example}
\begin{example}[\textbf{projection operator}]
\textbf{}\label{expansiveness_projection_operator}We remind the
reader that, a set is called \textbf{proximinal} if every point has
at least one projection onto it, whereas it is called a \textbf{Chebyshev}
set if every point has exactly one projection onto it. A nonempty
subset in $\mathbf{R}^{n}$ is Chebyshev if and only if it is closed
and convex \cite[Remark 3.15]{Bauschke2017}, and the projection operator
onto such a set is single-valued and firmly nonexpansive on $\mathbf{R}^{n}$
\cite[Proposition 4.8]{Bauschke2017}. However, for the constraint
set $\mathcal{C}$ in \eqref{eq:original_problem}, which is possibly
nonconvex, the projection operator $\tilde{\mathbf{\Pi}}_{\mathcal{C}}$
is not, in general, nonexpansive, hence not firmly nonexpansive. For
example, consider the set $\{0,1\}$; the projections of $0.4$ and
$0.6$ onto this set are 0 and 1, respectively, so $|0.6-0.4|=0.2<1$,
which violates the definition of nonexpansiveness. In such a case,
$2\tilde{\mathbf{\Pi}}_{\mathcal{C}}-I_{n}$ is also not nonexpansive,
because an operator is firmly nonexpansive if and only if its reflection
operator is nonexpansive \cite[Proposition 4.2]{Bauschke2017}. 

We now introduce the following definitions to (i) characterize an
operator that is not necessarily nonexpansive (\emph{e.g.,} $\tilde{\mathbf{\Pi}}_{\mathcal{C}}$
and $2\tilde{\mathbf{\Pi}}_{\mathcal{C}}-I_{n}$) and (ii) measure
the deviation of such an operator from being nonexpansive.
\end{example}

\paragraph*{Expansiveness of an operator\label{subsec:Expansiveness-of-an-op}}

Let $T:\mathcal{D}\to\mathbf{R}^{n}$ be a single-valued operator.
The\textbf{ expansiveness} of $T$ at $x,y$ in $\mathcal{D}$, denoted
by $\varepsilon_{xy}^{(T)}$, is defined as
\[
\varepsilon_{xy}^{(T)}=\begin{cases}
\|T(x)-T(y)\|-\|x-y\|,\\
\qquad\qquad\qquad\qquad\textrm{if }\|x-y\|<\|T(x)-T(y)\|\\
0,\qquad\qquad\qquad\quad{\rm else}.
\end{cases}
\]
where it is nonnegative and symmetric, \emph{i.e.,} $\varepsilon_{xy}^{(T)}=\varepsilon_{yx}^{(T)}\geq0$.
It follows that for every $x,y$ in $\mathcal{D}$,
\begin{equation}
\|T(x)-T(y)\|\leq\|x-y\|+\varepsilon_{xy}^{(T)}.\label{eq:N_x_y-1}
\end{equation}
{\footnotesize{}}Furthermore, define, \textbf{squared expansiveness}
of $T$ at $x,y$ in $\mathcal{D}$ as 
\[
\sigma_{xy}^{(T)}=\begin{cases}
\sqrt{\|T(x)-T(y)\|^{2}-\|x-y\|^{2}},\\
\qquad\qquad\qquad\qquad\textrm{if }\|x-y\|<\|T(x)-T(y)\|\\
0,\qquad\qquad\qquad\quad{\rm else}.
\end{cases}
\]
 Clearly, $\sigma_{xy}^{(T)}$ can be defined through $\varepsilon{}_{xy}^{(T)}$
as
\[
\sigma_{xy}^{(T)}=\sqrt{\varepsilon{}_{xy}^{(T)}}\sqrt{\|T(x)-T(y)\|+\|x-y\|}.
\]
It follows that for every $x,y$ in $\mathcal{D}$, 
\begin{equation}
\|T(x)-T(y)\|^{2}\leq\|x-y\|^{2}+\big(\sigma_{xy}^{(T)}\big)^{2}.\label{eq:M_x_y-1}
\end{equation}

\begin{rem}[\textbf{further characterization of nonexpansive operators}]
\label{(further-characterization-of}An operator $T$ is nonexpansive
on $\mathbf{R}^{n}$ if and only if $\varepsilon_{xy}^{(T)}=\sigma_{xy}^{(T)}=0$
for every $x,y$ in $\mathbf{R}^{n}$. On the other hand, an operator
$T$ is not nonexpansive if and only if there exist $x,y$ in its
domain such that $\varepsilon_{xy}^{(T)}$ is positive. Thus, $\varepsilon_{xy}^{(T)}$
and $\sigma_{xy}^{(T)}$ measure the deviation of $T$ from being
a nonexpansive operator at $x,y$.
\end{rem}

\subsection{Resolvent and reflected resolvent of an operator\label{subsec:Resolvent-and-reflected}}

Let $T:\mathbf{R}^{n}\rightrightarrows\mathbf{R}^{n}$ be a set-valued
operator and let $\gamma>0$. The \textbf{resolvent} of $T$, denoted
by $J_{\gamma T}$, is defined as $J_{\gamma T}=(I_{n}+\gamma T)^{-1}$,
and its \textbf{reflected resolvent,} denoted by $R_{\gamma T}$,
is defined as $R_{\gamma T}=2J_{\gamma T}-I_{n}$. The proximal operator
of a function is intimately connected to the resolvent of that function's
subdifferential operator as follows.
\begin{lem}[\textbf{resolvent characterization of proximal operator}]
\label{lem:singtonness}Let $g:\mathbf{R}^{n}\to\overline{\mathbf{R}}$
be proper, let $x\in\mathbf{R}^{n}$, and let $\gamma>0$. Then, both
$\mathbf{prox}_{\gamma g}$ and $J_{\gamma\partial g}$ are set-valued,
and $\mathbf{prox}_{\gamma g}(x)\subseteq J_{\gamma\partial g}(x)$.
Moreover, if $g$ is CPC, then both $\mathbf{prox}_{\gamma g}$ and
$J_{\gamma\partial g}$ are single-valued, firmly nonexpansive and
continuous on $\mathbf{R}^{n}$, and $\mathbf{prox}_{\gamma g}(x)=J_{\gamma\partial g}(x)$
. 
\end{lem}
\begin{IEEEproof}
When $g$ is proper, the claim follows from \cite[Example 10.2]{Rockafellar2009}.
When, $g$ is CPC, the claim follows from \cite[Proposition 12.27]{Bauschke2017},
\cite[pages 59-60]{Bauschke2017}, and \cite[Example 23.3]{Bauschke2017}. 
\end{IEEEproof}
The following corollary applies Lemma \ref{lem:singtonness} to the
constraint set $\mathcal{C}$ in \eqref{eq:original_problem}.
\begin{cor}[\textbf{resolvent characterization of projection}]
\emph{\label{lem:Proximal-operator-of}} For the constraint set $\mathcal{C}$
in \eqref{eq:original_problem}, $\tilde{\mathbf{\Pi}}_{\mathcal{C}}\left(x\right)\subseteq\mathbf{prox}_{\gamma\delta_{\mathcal{C}}}\left(x\right)=\mathbf{\Pi}_{\mathcal{C}}\left(x\right)\subseteq J_{\gamma\partial\delta_{C}}\left(x\right)$
for every $x\in\mathbf{R}^{n}$. For a convex set, all these operators
are single-valued, firmly nonexpansive, and equal to each other.
\end{cor}
\begin{IEEEproof}
Follows directly from Lemma \ref{lem:singtonness} and the definitions
of the proximal operator and the projection operator. 
\end{IEEEproof}

\section{Convergence of \eqref{eq:DR_splitting_orig}\label{sec:Convergence-analysis-for-NCDR}}

This section is organized as follows. First, in \S\ref{subsec:Supporting-lemmas-on},
we present some supporting lemmas on convergence of sequences. Then,
in \S\ref{subsec:Nonconvex-Douglas-Rachford}, we describe three
interrelated operators to develop the machinery for the convergence
analysis of \eqref{eq:DR_splitting_orig}, and in \S\ref{subsec:Characterization-of-global},
we characterize global minimizers of \eqref{eq:original_problem}
using these operators. In \S\ref{subsec:Main-convergence-result},
we present our main convergence result.

\subsection{Supporting lemmas on sequences\label{subsec:Supporting-lemmas-on}}

In this subsection, we present some supporting lemmas on sequences
to be used later; the first three results concern convergence of a
sequence of scalars, and the fourth result is about convergence of
a sequence of vectors in a compact set. 

First, we briefly review the definitions and basic properties of limit
inferior and limit superior of a sequence.\textbf{ Limit inferior}
and \textbf{limit superior} of a scalar sequence $\left(\alpha_{n}\right)_{n\in\mathbf{N}}$
are defined as
\begin{align*}
\underset{n\to\infty}{\underline{\lim}}\alpha_{n} & =\lim_{n\to\infty}\big(\inf_{m\geq n}\alpha_{m}\big),\textrm{and}\\
\underset{n\to\infty}{\overline{\lim}}\alpha_{n} & =\lim_{n\to\infty}\big(\sup_{m\geq n}\alpha_{m}\big),
\end{align*}
respectively, where they can be extended real-valued. For a bounded
sequence, both $\underline{\lim}_{n\to\infty}\alpha_{n}$ and $\overline{\lim}_{n\to\infty}\alpha_{n}$
exist, and they are finite. Clearly, $\underline{\lim}_{n\to\infty}\alpha_{n}\leq\overline{\lim}_{n\to\infty}\alpha_{n}$.
The sequence converges if and only if $\underline{\lim}_{n\to\infty}\alpha_{n}=\overline{\lim}_{n\to\infty}\alpha_{n}=\lim_{n\to\infty}\alpha_{n}\in\mathbf{R}.$
Furthermore, limit inferior satisfies \textbf{superadditivity}, \emph{i.e.,}
for every two sequences of real numbers, $\left(\alpha_{n}\right)_{n\in\mathbf{N}},\left(\beta_{n}\right)_{n\in\mathbf{N}}$
we have $\underline{\lim}_{n\to\infty}\left(\alpha_{n}+\beta_{n}\right)\geq\underline{\lim}_{n\to\infty}\alpha_{n}+\underline{\lim}_{n\to\infty}\beta_{n}.$ 
\begin{lem}[\textbf{limit of a nonnegative scalar sequence}]
\label{lem:(a-lemma-on} Let $\left(\alpha_{n}\right)_{n\in\mathbf{N}}$
be a sequence of nonnegative scalars such that $\sum_{n\in\mathbf{N}}\alpha_{n}$
is bounded above. Then, $\lim_{n\to\infty}\alpha_{n}=0$.
\end{lem}
\begin{IEEEproof}
Directly follows from \cite[Proposition 3.2.1]{Davidson09} and \cite[Theorem 3.1.4]{Davidson09}.
\end{IEEEproof}
\begin{lem}[{\textbf{convergence of a nonnegative scalar sequence }\cite[page 44, Lemma 2]{Polyak1987}}]
\emph{\label{lem:Polyak Lemma}} Let $\left(u_{n}\right)_{n\in\mathbf{N}},\left(\alpha_{n}\right)_{n\in\mathbf{N}}$,
and $\left(\beta_{n}\right)_{n\in\mathbf{N}}$ be sequences of nonnegative
scalars such that for every $n\in\mathbf{N}$, we have $u_{n+1}\leq(1+\alpha_{n})u_{n}+\beta_{n},$
$\sum_{n\in\mathbf{N}}\alpha_{n}<\infty,$ and $\sum_{n\in\mathbf{N}}\beta_{n}<\infty.$
Then, there is a nonnegative scalar $u$ such that $u_{n}$ converges
to $u$. 
\end{lem}

\begin{lem}[{\textbf{limit inferior of addition of two sequences }\cite[Proposition 2.3]{bworld}}]
\label{lem:-Let-df} Let $\left(\alpha_{n}\right)_{n\in\mathbf{N}}$
and $\left(\beta_{n}\right)_{n\in\mathbf{N}}$ be two bounded scalar
sequences. If $\lim_{n\to\infty}\alpha_{n}=\alpha$, then $\mathop{\underline{\mathrm{lim}}}_{n\to\infty}\left(\alpha_{n}+\beta_{n}\right)=\alpha+\mathop{\underline{\mathrm{lim}}}_{n\to\infty}\beta_{n}.$
\end{lem}
Now we record a result about convergence of a sequence of vectors
in a compact set. We remind the reader that, a set is \textbf{connected}
if it is not the union of two disjoint nonempty closed sets. A compact
and connected set is called a \textbf{continuum}. Moreover, a set
is called a \textbf{nontrivial continuum}, if it is a continuum, and
it does not reduce to $\emptyset$ or a singleton \cite{Combettes1990}.
Finally, $x$ is a \textbf{cluster point} of a sequence $\left(x_{n}\right)_{n\in\mathbf{N}}$
if the sequence has a subsequence that converges to $x$. 
\begin{lem}[{\textbf{convergence of a sequence of vectors in a compact set }\cite[Theorem 4.2]{Combettes1990}}]
\emph{} \label{lem:Combettes_Proof} Let $\left(x_{n}\right)_{n\in\mathbf{N}}$
be a sequence of vectors in a compact set $\mathcal{S}\subseteq\mathbf{R}^{n}$
such that $\|x_{n+1}-x_{n}\|$ converges to zero. Then, either $\left(x_{n}\right)_{n\in\mathbf{N}}$
converges to a point in $\mathcal{S}$, or its set of cluster points
is a nontrivial continuum in $\mathcal{S}$.
\end{lem}

\subsection{Nonconvex Douglas-Rachford, Cayley and Peaceman-Rachford operators\label{subsec:Nonconvex-Douglas-Rachford}}

To facilitate our convergence analysis, we define the following operators
for \eqref{eq:original_problem}.
\begin{itemize}
\item The \textbf{nonconvex Douglas-Rachford operator} with parameter $\gamma>0$,
denoted by $\tilde{T}$, is defined as 
\begin{equation}
\tilde{T}=\tilde{\mathbf{\Pi}}_{\mathcal{\mathcal{C}}}\left(2\mathbf{prox}_{\gamma f}-I_{n}\right)+I_{n}-\mathbf{prox}_{\gamma f}.\label{eq:NC-DRS-operator}
\end{equation}
\item The \textbf{nonconvex Cayley operator} of $\tilde{T}$ (also known
as the reflection operator of $\tilde{T}$) with parameter $\gamma>0$,
denoted by $\tilde{R}$, is defined as 
\begin{equation}
\tilde{R}=2\tilde{T}-I_{n}.\label{eq:NC-Cayley-operator}
\end{equation}
\item The \textbf{nonconvex Peaceman-Rachford operator} with parameter $\gamma>0$,
denoted by $\tilde{S}$, is defined as 
\begin{equation}
\tilde{S}=(2\mathbf{\tilde{\Pi}}_{\mathcal{C}}-I_{n})(2\mathbf{prox}_{\gamma f}-I_{n}).\label{eq:NC-PR-operator}
\end{equation}
\end{itemize}
\begin{rem}[\textbf{nonconvex Peaceman-Rachford operator $\tilde{S}$ is not nonexpansive}]
\label{tilde_S_nonexpansive} Note that $\tilde{S}$ is a composition
of $2\mathbf{\tilde{\Pi}}_{\mathcal{C}}-I_{n}$ and $2\mathbf{prox}_{\gamma f}-I_{n}$,
where the latter is nonexpansive (see Example \ref{(proximal-operator)-The}),
but the former is not nonexpansive in general (see Example \ref{expansiveness_projection_operator}).
Hence $\tilde{S}$ is not a nonexpansive operator in general.
\end{rem}
These operators allow us to write \eqref{eq:DR_splitting_orig} in
the following compact form:
\begin{align}
z_{n+1}=\tilde{T}z_{n} & =\frac{1}{2}\tilde{R}z_{n}+\frac{1}{2}z_{n}.\tag{Compact-NC-DRS}\label{eq:convenient-form-1}
\end{align}

The following lemma will be used later to characterize global minimizers
of \eqref{eq:original_problem}. 
\begin{lem}[\textbf{characterization of nonconvex Peaceman-Rachford operator}]
\label{tilde_S_in_R_R} For \eqref{eq:original_problem}, let $\tilde{S}$
be the nonconvex Peaceman-Rachford operator with parameter $\gamma>0$
defined in \eqref{eq:NC-PR-operator}. Then, $\tilde{S}\left(x\right)\subseteq R_{\gamma\partial\delta_{\mathcal{C}}}R_{\gamma\partial f}\left(x\right)$
for every $x\in\mathbf{R}^{n}$. 
\end{lem}
\begin{IEEEproof}
As $f$ is CPC, we have 
\begin{equation}
R_{\gamma\partial f}=2\mathbf{prox}_{\gamma f}-I_{n},\label{eq:asgard}
\end{equation}
using Lemma \ref{lem:singtonness} and the definition of the reflected
resolvent in $\S$\ref{subsec:Resolvent-and-reflected}. Now for every
$x\in\mathbf{R}^{n}$,
\begin{align}
\big(2\tilde{\mathbf{\Pi}}_{\mathcal{C}}-I_{n}\big)\left(x\right) & =2\tilde{\mathbf{\Pi}}_{\mathcal{C}}\left(x\right)-x\nonumber \\
 & \overset{a)}{\subseteq}2\mathbf{\Pi}_{\mathcal{C}}\left(x\right)-x\nonumber \\
 & \overset{b)}{\subseteq}2J_{\gamma\delta_{\mathcal{C}}}(x)-x\nonumber \\
 & =\big(2J_{\gamma\delta_{\mathcal{C}}}-I_{n}\big)x\nonumber \\
 & \overset{c)}{=}R_{\gamma\delta_{\mathcal{C}}}\left(x\right),\label{eq:midgard}
\end{align}
where $a)$ follows from $\tilde{\mathbf{\Pi}}_{\mathcal{C}}\left(x\right)\subseteq\mathbf{\Pi}_{\mathcal{C}}\left(x\right)$
for every $x$ in $\mathbf{R}^{n}$ (Example \ref{exa:(Projection-operator)}),
$b)$ follows from Corollary \ref{lem:Proximal-operator-of}, and
$c)$ follows from the definition of reflected resolvent in \S \ref{subsec:Resolvent-and-reflected}.
Thus, for every $x\in\mathbf{R}^{n}$,
\begin{align*}
\tilde{S}\left(x\right) & =\big(2\mathbf{\tilde{\Pi}}_{\mathcal{C}}-I_{n}\big)\big(2\mathbf{prox}_{\gamma f}-I_{n}\big)\left(x\right)\\
 & \overset{a)}{=}\big(2\mathbf{\tilde{\Pi}}_{\mathcal{C}}-I_{n}\big)R_{\gamma\partial f}\left(x\right)\\
 & \overset{b)}{\subseteq}R_{\gamma\partial\delta_{\mathcal{C}}}R_{\gamma\partial f}\left(x\right),
\end{align*}
 where $a)$ and $b)$ use \eqref{eq:asgard} and \eqref{eq:midgard},
respectively.
\end{IEEEproof}
\begin{prop}[\textbf{relationship between $\tilde{T}$, $\tilde{R}$, and $\tilde{S}$}]
\label{prop:(Relationship-between-the} For \eqref{eq:original_problem},
let $\tilde{T}$, $\tilde{R}$, and $\tilde{S}$ be the operators
with parameter $\gamma>0$ defined in \eqref{eq:NC-DRS-operator},
\eqref{eq:NC-Cayley-operator}, and \eqref{eq:NC-PR-operator}, respectively.
Then, 

\emph{(i)} the operators $\tilde{R}$ and $\tilde{S}$ are equal,
i.e., $\tilde{R}\left(x\right)=\tilde{S}\left(x\right)$ for every
$x\in\mathbf{R}^{n}$, and 

\emph{(ii)} the fixed point sets of $\tilde{T}$, $\tilde{R}$, and
$\tilde{S}$ are equal, i.e., $\mathop{{\bf fix}}\tilde{R}=\mathop{{\bf fix}}\tilde{S}=\mathop{{\bf fix}}\tilde{T}$.
\end{prop}
\begin{IEEEproof}
(i): For every $x\in\mathbf{R}^{n},$ 
\begin{align}
\tilde{T}\left(x\right) & =\big(\mathbf{\tilde{\Pi}}_{\mathcal{C}}\big(2\mathbf{prox}_{\gamma f}-I_{n}\big)+I_{n}-\mathbf{prox}_{\gamma f}\big)\left(x\right)\nonumber \\
 & =\mathbf{\tilde{\Pi}}_{\mathcal{C}}\big(2\mathbf{prox}_{\gamma f}-I_{n}\big)\left(x\right)+x-\mathbf{prox}_{\gamma f}\left(x\right)\nonumber \\
 & =\mathbf{\tilde{\Pi}}_{\mathcal{C}}\big(2\mathbf{prox}_{\gamma f}\left(x\right)-x\big)+x-\mathbf{prox}_{\gamma f}\left(x\right).\label{eq:obs_douglas_rachford}
\end{align}
Furthermore, for every $x\in\mathbf{R}^{n}$, 
\begin{align}
\tilde{R}\left(x\right) & =\big(2\tilde{T}-I_{n}\big)\left(x\right)\nonumber \\
 & =2\tilde{T}\left(x\right)-x\nonumber \\
 & \overset{a)}{=}2\mathbf{\tilde{\Pi}}_{\mathcal{C}}\left(2\mathbf{prox}_{\gamma f}\left(x\right)-x\right)+2x-2\mathbf{prox}_{\gamma f}\left(x\right)-x\nonumber \\
 & =2\mathbf{\tilde{\Pi}}_{\mathcal{C}}\left(2\mathbf{prox}_{\gamma f}\left(x\right)-x\right)+x-2\mathbf{prox}_{\gamma f}\left(x\right),\label{eq:obs_douglas_rachford_2}
\end{align}
 where $a)$ uses \eqref{eq:obs_douglas_rachford}. Hence, for every
$x\in\mathbf{R}^{n}$
\begin{align*}
\tilde{S}\left(x\right) & =\big(2\mathbf{\tilde{\Pi}}_{\mathcal{C}}-I_{n}\big)\big(2\mathbf{prox}_{\gamma f}-I_{n}\big)\left(x\right)\\
 & =\big(2\mathbf{\tilde{\Pi}}_{\mathcal{C}}-I_{n}\big)\underbrace{\big(2\mathbf{prox}_{\gamma f}\left(x\right)-x\big)}_{=y\textrm{ (let)}}\\
 & =2\mathbf{\tilde{\Pi}}_{\mathcal{C}}\left(y\right)-y\\
 & =2\mathbf{\tilde{\Pi}}_{\mathcal{C}}\big(2\mathbf{prox}_{\gamma f}\left(x\right)-x\big)-2\mathbf{prox}_{\gamma f}\left(x\right)+x\\
 & \overset{a)}{=}\tilde{R}x,
\end{align*}
 where $a)$ uses \eqref{eq:obs_douglas_rachford_2}. 

(ii): In (i), $\tilde{R}=\tilde{S}$ implying $\mathop{{\bf fix}}\tilde{R}=\mathop{{\bf fix}}\tilde{S}$.
Now $x\in\mathop{{\bf fix}}\tilde{T}\Leftrightarrow\tilde{T}\left(x\right)=x\Leftrightarrow2\tilde{T}\left(x\right)=2x\Leftrightarrow2\tilde{T}\left(x\right)-x=x\Leftrightarrow\left(2\tilde{T}-I_{n}\right)\left(x\right)=x\Leftrightarrow x=\mathop{{\bf fix}}\tilde{R}.$
So $\mathop{{\bf fix}}\tilde{T}=\mathop{{\bf fix}}\tilde{R}=\mathop{{\bf fix}}\tilde{S}$. 
\end{IEEEproof}
\begin{rem}[\textbf{nonconvex Cayley operator $\tilde{R}$ is not nonexpansive}]
\label{Cayley_is_not_nonexpansive} From Remark \ref{tilde_S_nonexpansive}
and Proposition \ref{prop:(Relationship-between-the}, it follows
that $\tilde{R}$ is not a nonexpansive operator in general. This
plays an important role in our convergence analysis; in particular,
the sufficient conditions for convergence of \eqref{eq:DR_splitting_orig}
are dictated by the squared expansiveness of $\tilde{R}$ over the
iterates of \eqref{eq:DR_splitting_orig}.
\end{rem}

\subsection{Characterization of global minimizers\label{subsec:Characterization-of-global}}

Minimizers of \eqref{eq:original_problem} are characterized via the
nonconvex Douglas-Rachford operator as follows. 
\begin{thm}[\textbf{global minimizers of \eqref{eq:original_problem}}]
\emph{\label{prop:minimizer_charecterization}} For \eqref{eq:original_problem},
let $\tilde{T}$ be the nonconvex Douglas-Rachford operator with parameter
$\gamma>0$ defined in \eqref{eq:NC-DRS-operator}. Then,

\emph{(i)} sum of the functions $f+\delta_{\mathcal{C}}$ is proper,
$\partial f\left(x\right)+\partial\delta_{\mathcal{C}}\left(x\right)\subseteq\partial\left(f+\delta_{\mathcal{C}}\right)\left(x\right)$
for every $x\in\mathbf{R}^{n}$, and both $\partial f+\partial\delta_{\mathcal{C}}$
and $\partial\left(f+\delta_{\mathcal{C}}\right)$ are monotone operators, 

\emph{(ii) $\mathop{{\bf zer}}\big(\partial f+\partial\delta_{\mathcal{C}}\big)=\mathbf{prox}_{\gamma f}\big(\mathop{{\bf fix}}(R_{\gamma\partial\delta_{\mathcal{C}}}R_{\gamma\partial f})\big)$},
and

\emph{(iii)} if $\mathop{{\bf fix}}\tilde{T}\neq\emptyset$, then
$\mathbf{prox}_{\gamma f}\big(\mathop{{\bf fix}}\tilde{T}\big)\subseteq\mathop{{\rm argmin}}\left(f+\delta_{\mathcal{C}}\right)$.
\end{thm}
\begin{IEEEproof}
(i): The indicator function of a closed set is closed \cite[Example 1.25]{Bauschke2017},
and the indicator function of a nonempty set is proper \cite[pages 6-7]{Rockafellar2009}.
Hence, $\delta_{\mathcal{C}}$ is closed and proper. Also, we have
$\mathop{{\bf dom}}f\cap\mathop{{\bf dom}}\delta_{\mathcal{C}}\neq\emptyset$,
otherwise \eqref{eq:original_problem} is infeasible. So, using Lemma
\ref{lem:subdifferential-monotone}, the function $f+\delta_{\mathcal{C}}$
is proper, $\partial f\left(x\right)+\partial\delta_{\mathcal{C}}\left(x\right)\subseteq\partial\left(f+\delta_{\mathcal{C}}\right)\left(x\right)$
for every $x\in\mathbf{R}^{n}$, and both $\partial f+\partial\delta_{\mathcal{C}}$
and $\partial\left(f+\delta_{\mathcal{C}}\right)$ are monotone operators.

(ii): This proof is based on \cite[Proposition 25.1 (ii)]{Bauschke2017}.
For every $\gamma>0$, we have
\begin{flalign}
 & x\in\mathop{{\bf zer}}\big(\partial f+\partial\delta_{\mathcal{C}}\big)\nonumber \\
\Leftrightarrow & (\exists y\in\mathbf{R}^{n})\;x-y\in\gamma\partial\delta_{\mathcal{C}}(x)\textrm{ and }y-x\in\gamma\partial f(x)\nonumber \\
\Leftrightarrow & \left(\exists y\in\mathbf{R}^{n}\right)\;2x-y\in(I_{n}+\gamma\partial\delta_{\mathcal{C}})(x)\textrm{ and }\nonumber \\
 & \mbox{\qquad\qquad\quad}y\in(I_{n}+\gamma\partial f)(x)\nonumber \\
\Leftrightarrow & \left(\exists y\in\mathbf{R}^{n}\right)\;\underbrace{(I_{n}+\gamma\partial\delta_{\mathcal{C}})^{-1}}_{=J_{\gamma\partial\delta_{\mathcal{C}}}}(2x-y)\ni x\textrm{ and }\nonumber \\
 & \mbox{\qquad\qquad\quad}\underbrace{(I_{n}+\gamma\partial f)^{-1}}_{=J_{\gamma\partial f}}(y)\ni x\nonumber \\
\overset{a)}{\Leftrightarrow} & \left(\exists y\in\mathbf{R}^{n}\right)\;x\in J_{\gamma\partial\delta_{\mathcal{C}}}(2x-y)\textrm{ and }x=J_{\gamma\partial f}(y),\nonumber \\
\overset{b)}{\Leftrightarrow} & \left(\exists y\in\mathbf{R}^{n}\right)\;x\in J_{\gamma\partial\delta_{\mathcal{C}}}R_{\gamma\partial f}(y)\textrm{ and }x=J_{\gamma\partial f}(y)\label{eq:shapiro}
\end{flalign}
 where $a)$ uses the facts that $J_{\gamma\partial f}$ is a single-valued
operator (from Lemma \ref{lem:singtonness}), and $J_{\gamma\partial\delta_{\mathcal{C}}}$
is a set-valued operator (from Corollary \ref{lem:Proximal-operator-of}),
and $b)$ uses the observation that $x=J_{\gamma\partial f}(y)$ can
be expressed as
\[
x=J_{\gamma\partial f}(y)\Leftrightarrow2x-y=\big(2J_{\gamma\partial f}-I_{n}\big)y=R_{\gamma\partial f}(y).
\]
Also, using the last expression, we can write the first term of \eqref{eq:shapiro}
as
\begin{align}
 & J_{\gamma\partial\delta_{\mathcal{C}}}R_{\gamma\partial f}(y)\ni x\nonumber \\
\Leftrightarrow & 2J_{\gamma\partial\delta_{\mathcal{C}}}R_{\gamma\partial f}(y)-y\ni2x-y=R_{\gamma\partial f}(y)\nonumber \\
\Leftrightarrow & y\in2J_{\gamma\partial\delta_{\mathcal{C}}}R_{\gamma\partial f}(y)-R_{\gamma\partial f}(y)\nonumber \\
 & \;=\big(2J_{\gamma\partial\delta_{\mathcal{C}}}-I_{n}\big)\big(R_{\gamma\partial f}(y)\big)\nonumber \\
 & \;=R_{\gamma\partial\delta_{\mathcal{C}}}R_{\gamma\partial f}(y)\nonumber \\
\Leftrightarrow & y\in\mathop{{\bf fix}}\big(R_{\gamma\partial\delta_{\mathcal{C}}}R_{\gamma\partial f}\big).\label{eq:harris}
\end{align}
Using \eqref{eq:shapiro}, \eqref{eq:harris}, and $J_{\gamma\partial f}=\mathbf{prox}_{\gamma f}$
(from Lemma \ref{lem:singtonness}) we have
\begin{align*}
 & x\in\mathop{{\bf zer}}\big(\partial f+\partial\delta_{\mathcal{C}}\big)\\
\Leftrightarrow & \left(\exists y\in\mathbf{R}^{n}\right)\;y\in\mathop{{\bf fix}}\big(R_{\gamma\partial\delta_{\mathcal{C}}}R_{\gamma\partial f}\big)\textrm{ and }x=\mathbf{prox}_{\gamma f}(y)\\
\Leftrightarrow & x\in\mathbf{prox}_{\gamma f}\big(\mathop{{\bf fix}}(R_{\gamma\partial\delta_{\mathcal{C}}}R_{\gamma\partial f})\big).
\end{align*}
Thus, \emph{$\mathop{{\bf zer}}\big(\partial f+\partial\delta_{\mathcal{C}}\big)=\mathbf{prox}_{\gamma f}\big(\mathop{{\bf fix}}(R_{\gamma\partial\delta_{\mathcal{C}}}R_{\gamma\partial f})\big).$}

(iii): We have
\begin{align*}
 & x\in\mathop{{\bf zer}}\left(\partial f+\partial\delta_{\mathcal{C}}\right)\\
\Leftrightarrow & 0\in\partial f\left(x\right)+\partial\delta_{\mathcal{C}}\left(x\right)\overset{a)}{\subseteq}\partial\left(f+\delta_{\mathcal{C}}\right)\left(x\right)\\
\Rightarrow & x\in\mathop{{\bf zer}}\partial\left(f+\delta_{\mathcal{C}}\right),
\end{align*}
where $a)$ uses $\partial f\left(x\right)+\partial\delta_{\mathcal{C}}\left(x\right)\subseteq\partial\left(f+\delta_{\mathcal{C}}\right)\left(x\right)$
proven in (i). So, $\mathop{{\bf zer}}\left(\partial f+\partial\delta_{\mathcal{C}}\right)\subseteq\mathop{{\bf zer}}\left(\partial\left(f+\delta_{\mathcal{C}}\right)\right)$.
Combining the last statement with $\mathop{{\bf zer}}\left(\partial\left(f+\delta_{\mathcal{C}}\right)\right)=\mathop{{\rm argmin}}\left(f+\delta_{\mathcal{C}}\right)$
(from Theorem \ref{thm:(Fermat's-rule)}) and (ii), we have
\begin{align}
 & \mathop{{\bf zer}}\big(\partial f+\partial\delta_{\mathcal{C}}\big)\nonumber \\
= & \mathbf{prox}_{\gamma f}\Big(\mathop{{\bf fix}}\big(R_{\gamma\partial\delta_{\mathcal{C}}}R_{\gamma\partial f}\big)\Big)\\
\subseteq & \mathop{{\rm argmin}}\big(f+\delta_{\mathcal{C}}\big).\label{eq:nifelhaim}
\end{align}
Recall from Lemma \ref{lem:singtonness} that $\mathbf{prox}_{\gamma f}$
is a single-valued operator. Thus, 
\begin{align*}
\mathbf{prox}_{\gamma f}\big(\mathop{{\bf fix}}\tilde{S}\big) & =\bigcup_{x\in\mathop{{\bf fix}}\tilde{S}}\mathbf{prox}_{\gamma f}\left(x\right)\\
 & =\bigcup_{x:x=\tilde{S}(x)}\mathbf{prox}_{\gamma f}\left(x\right)\\
 & \overset{a)}{\subseteq}\bigcup_{x:x\in R_{\gamma\partial\delta_{\mathcal{C}}}R_{\gamma\partial f}\left(x\right)}\mathbf{prox}_{\gamma f}\left(x\right)\\
 & =\bigcup_{x:x\in\mathop{{\bf fix}}R_{\gamma\partial\delta_{\mathcal{C}}}R_{\gamma\partial f}}\mathbf{prox}_{\gamma f}\left(x\right)\\
 & =\mathbf{prox}_{\gamma f}\big(\mathop{{\bf fix}}R_{\gamma\partial\delta_{\mathcal{C}}}R_{\gamma\partial f}\big),
\end{align*}
where $a)$ uses $\tilde{S}\left(x\right)\subseteq R_{\gamma\partial\delta_{\mathcal{C}}}R_{\gamma\partial f}\left(x\right)$
for every $x\in\mathbf{R}^{n}$ (from Lemma \ref{tilde_S_in_R_R}).
But, $\mathop{{\bf fix}}\tilde{S}=\mathop{{\bf fix}}\tilde{T}$ from
Proposition \ref{prop:(Relationship-between-the}. So, 
\begin{align*}
\mathbf{prox}_{\gamma f}\big(\mathop{{\bf fix}}\tilde{S}\big) & =\mathbf{prox}_{\gamma f}\big(\mathop{{\bf fix}}\tilde{T}\big)\\
 & \subseteq\mathbf{prox}_{\gamma f}\big(\mathop{{\bf fix}}R_{\gamma\partial\delta_{\mathcal{C}}}R_{\gamma\partial f}\big).
\end{align*}
Combining the last equation with \eqref{eq:nifelhaim}, we have
\begin{align*}
\mathbf{prox}_{\gamma f}\big(\mathop{{\bf fix}}\tilde{T}\big) & \subseteq\mathbf{prox}_{\gamma f}\big(\mathop{{\bf fix}}R_{\gamma\partial\delta_{\mathcal{C}}}R_{\gamma\partial f}\big)\\
 & \subseteq\mathop{{\rm argmin}}\big(f+\delta_{\mathcal{C}}\big).
\end{align*}
\end{IEEEproof}
\begin{rem}[\textbf{nonemptiness of $\mathop{{\bf zer}}\big(\partial f+\partial\delta_{\mathcal{C}}\big)$}]
 A necessary condition for nonemptiness of $\mathop{{\bf fix}}\tilde{T}$
is nonemptiness of $\mathop{{\bf zer}}\big(\partial f+\partial\delta_{\mathcal{C}}\big)$.
This necessary condition $\mathop{{\bf zer}}\big(\partial f+\partial\delta_{\mathcal{C}}\big)\neq\emptyset$
is stronger than the existence of a minimizer, because, even in a
convex setup, $\mathop{{\bf zer}}\big(\partial f+\partial\delta_{\mathcal{C}}\big)\neq\mathop{{\bf zer}}\left(\partial\left(f+\delta_{\mathcal{C}}\right)\right)$,
in general \cite[Remark 16.7]{Bauschke2017}. Nevertheless, we will
assume that $\mathop{{\bf zer}}\big(\partial f+\partial\delta_{\mathcal{C}}\big)\neq\emptyset$
for the rest of our development, as this seems to be a standard assumption
even in convex optimization literature \cite{Davis2016}.
\end{rem}

\subsection{Main convergence result\label{subsec:Main-convergence-result}}

We remind the reader that the nonconvex Cayley operator $\tilde{R}$
is not nonexpansive in general (Remark \ref{Cayley_is_not_nonexpansive}).
To characterize the deviation of $\tilde{R}$ from being a nonexpansive
operator, recalling $\S$\ref{subsec:Nonexpansive,-firmly-nonexpansiv},
we use expansiveness and squared expansiveness of $\tilde{R}$ at
each $x,y$ in $\mathbf{R}^{n}$, denoted by $\varepsilon_{xy}$ and
$\sigma_{xy}$, respectively; here we have dropped the superscript
$(\tilde{R})$ to reduce notational burden. So, from \eqref{eq:N_x_y-1}
and \eqref{eq:M_x_y-1}, for every $x,y$ in $\mathbf{R}^{n}$,
\begin{align}
\|\tilde{R}(x)-\tilde{R}(y)\| & \leq\|x-y\|+\varepsilon_{xy},\textrm{and}\label{eq:N_x_y}\\
\|\tilde{R}(x)-\tilde{R}(y)\|^{2} & \leq\|x-y\|^{2}+\sigma_{xy}^{2}.\label{eq:M_x_y}
\end{align}

Also, the closed ball with center $x\in\mathbf{R}^{n}$ and finite
radius $r>0$, denoted by $B(x;r)$, is defined as $B(x;r)=\{y\mid\|x-y\|\leq r\}$;
a closed ball in $\mathbf{R}^{n}$ with finite radius is compact \cite[\S 2.4]{Bauschke2017}.
Now we present our main convergence result. 
\begin{thm}[\textbf{main convergence result}]
\label{(main-convergence-result)} For \eqref{eq:original_problem},
let $(z_{n})_{n\in\mathbf{N}}$ be the sequence of vectors generated
by \eqref{eq:DR_splitting_orig}. Suppose that, for the chosen initial
point $z_{0}$, there exists a $z\in\mathop{{\bf fix}}\tilde{T}$,
such that $\sum_{n=0}^{\infty}\sigma_{z_{n}z}^{2}$ is bounded above,
and $\|z_{0}-z\|^{2}$ is finite.\textbf{} Define $r:=\sqrt{\|z_{0}-z\|^{2}+\frac{1}{2}\sum_{n=0}^{\infty}\sigma_{z_{n}z}^{2}}$.
Then, one of the following holds:

\emph{(i)} the sequence $(z_{n})_{n\in\mathbf{N}}$ converges to a
point $z^{\star}\in B(z;r)$. In this case, suppose also that $\mathop{\underline{\mathrm{lim}}}_{n\to\infty}\sigma_{z_{n}z^{\star}}^{2}=0$.
Then, $\mathbf{prox}_{\gamma f}\left(z^{\star}\right)$ is an optimal
solution of \eqref{eq:original_problem}, and the sequence $(x_{n})_{n\in\mathbf{N}}$
generated by \eqref{eq:DR_splitting_orig} converges to $\mathbf{prox}_{\gamma f}(z^{\star})$.

\emph{(ii)} the set of cluster points of $(z_{n})_{n\in\mathbf{N}}$
forms a nontrivial continuum in $B(z;r)$.
\end{thm}
\begin{IEEEproof}
\textbf{Step 1.} First, we show that the sequence $(z_{n})_{n\in\mathbf{N}}$
stays in the compact set $B(z;r)$. For every $n\in\mathbf{N}$, 
\begin{align}
\|z_{n+1}-z\|^{2} & \overset{a)}{=}\|z_{n}+\frac{1}{2}\big(\tilde{R}z_{n}-z_{n}\big)-z\|^{2}\nonumber \\
 & =\|\frac{1}{2}\big(z_{n}-z\big)+\frac{1}{2}\big(\tilde{R}z_{n}-z\big)\|^{2}\nonumber \\
 & \overset{b)}{=}\frac{1}{2}\|z_{n}-z\|^{2}+\frac{1}{2}\|\tilde{R}z_{n}-z\|^{2}\nonumber \\
 & \quad-\frac{1}{4}\|(z_{n}-z)-(\tilde{R}z_{n}-z)\|^{2}\nonumber \\
 & \overset{c)}{=}\frac{1}{2}\|z_{n}-z\|^{2}+\frac{1}{2}\|\tilde{R}z_{n}-\tilde{R}z\|^{2}\nonumber \\
 & \quad-\frac{1}{4}\|z_{n}-\tilde{R}z_{n}\|^{2}\nonumber \\
 & \overset{d)}{\leq}\frac{1}{2}\|z_{n}-z\|^{2}+\frac{1}{2}\|z_{n}-z\|^{2}\nonumber \\
 & \quad+\frac{1}{2}\sigma_{z_{n}z}^{2}-\frac{1}{4}\|z_{n}-\tilde{R}z_{n}\|^{2}\nonumber \\
 & \leq\|z_{n}-z\|^{2}-\frac{1}{4}\|z_{n}-\tilde{R}z_{n}\|^{2}+\frac{1}{2}\sigma_{z_{n}z}^{2}\label{eq:inequality_1-1-1}\\
 & \overset{e)}{\leq}\|z_{n}-z\|^{2}+\frac{1}{2}\sigma_{z_{n}z}^{2},\label{eq:aha_ki_anondo}
\end{align}
where $a)$ uses \eqref{eq:convenient-form-1}, $b)$ uses the identity
$\|\alpha x+(1-\alpha)y\|^{2}=\alpha\|x\|^{2}+(1-\alpha)\|y\|^{2}-\alpha(1-\alpha)\|x-y\|^{2}$
for every $x,y\in\mathbf{R}^{n}$ and every $\alpha\in\mathbf{R}$
\cite[Corollary 2.14]{Bauschke2017}, $c)$ uses $z\in\mathop{{\bf fix}}\tilde{T}$,
and $\mathop{{\bf fix}}\tilde{T}=\mathop{{\bf fix}}\tilde{R}$ (from
Proposition \ref{prop:(Relationship-between-the}(ii)), $d)$ uses
\eqref{eq:M_x_y}, and $e)$ is obtained by removing the nonpositive
term $-\frac{1}{4}\|z_{n}-\tilde{R}z_{n}\|^{2}$. From \eqref{eq:aha_ki_anondo},
we have
\begin{align}
\|z_{n}-z\|^{2} & \leq\|z_{n-1}-z\|^{2}+\frac{1}{2}\sigma_{z_{n-1}z}^{2}\nonumber \\
 & \leq\|z_{n-2}-z\|^{2}+\frac{1}{2}\sigma_{z_{n-2}z}^{2}+\frac{1}{2}\sigma_{z_{n-1}z}^{2}\nonumber \\
 & \leq\|z_{0}-z\|^{2}+\frac{1}{2}\sum_{i=0}^{n-1}\sigma_{z_{i}z}^{2}\nonumber \\
 & \leq\|z_{0}-z\|^{2}+\frac{1}{2}\sum_{i=0}^{\infty}\sigma_{z_{i}z}^{2},\label{eq:tyrs_temple}
\end{align}
where the final term is bounded, because $\sum_{i=0}^{\infty}\sigma_{z_{i}z}^{2}$
is bounded above. Hence, the sequence $(z_{n})_{n\in\mathbf{N}}$
stays in the compact set $B(z;r)$.

\textbf{Step 2. }Next, we show that $\lim_{n\to\infty}\|\tilde{R}z_{n}-z_{n}\|=0$.
From \eqref{eq:inequality_1-1-1},
\begin{align*}
\frac{1}{4}\|\tilde{R}z_{n}-z_{n}\|^{2} & \leq\big(\|z_{n}-z\|^{2}-\|z_{n+1}-z\|^{2}\big)+\frac{1}{2}\sigma_{z_{n}z}^{2}\\
\Rightarrow\frac{1}{4}\sum_{n=0}^{m}\|\tilde{R}z_{n}-z_{n}\|^{2} & \leq\sum_{n=0}^{m}\big(\|z_{n}-z\|^{2}-\|z_{n+1}-z\|^{2}\big)\\
 & \quad+\frac{1}{2}\sum_{n=0}^{m}\sigma_{z_{n}z}^{2}\\
 & \overset{a)}{=}\big(\|z_{0}-z\|^{2}-\|z_{m+1}-z\|^{2}\big)\\
 & \quad+\frac{1}{2}\sum_{n=0}^{m}\sigma_{z_{n}z}^{2}\\
 & \overset{b)}{\leq}\|z_{0}-z\|^{2}+\frac{1}{2}\sum_{n=0}^{m}\sigma_{z_{n}z}^{2},
\end{align*}
where $a)$ uses the telescopic sum, and $b)$ is obtained by removing
the negative term $\|z_{m+1}-z\|^{2}$. If $m\to\infty$, then the
right hand side of the last inequality is bounded above, because $\sum_{n=0}^{\infty}\sigma_{z_{n}z}^{2}$
is bounded above. Thus, $\sum_{n=0}^{\infty}\|\tilde{R}z_{n}-z_{n}\|^{2}$
is bounded above, and using Lemma \ref{lem:(a-lemma-on}, we have
$\lim_{n\to\infty}\|\tilde{R}z_{n}-z_{n}\|^{2}=0$, \emph{i.e.,} $\lim_{n\to\infty}\|\tilde{R}z_{n}-z_{n}\|=0$.

\textbf{Step 3. }We show that sequence $(z_{n})_{n\in\mathbf{N}}$
either converges to a point or its set of cluster points forms a nontrivial
continuum. In step 2, we have shown that, $\lim_{n\to\infty}\|\tilde{R}z_{n}-z_{n}\|=0.$
On the other hand, $\|\tilde{R}z_{n}-z_{n}\|=2\|z_{n+1}-z_{n}\|$
from \eqref{eq:convenient-form-1}, so $\lim_{n\to\infty}\|z_{n+1}-z_{n}\|=0$.
Thus, the sequence $\left(z_{n}\right)_{n\in\mathbf{N}}$ stays in
a compact set $B(z;r)$ and satisfies $\lim_{n\to\infty}\|z_{n+1}-z_{n}\|=0$.
So, due to Lemma \ref{lem:Combettes_Proof}, the sequence $(z_{n})_{n\in\mathbf{N}}$
either converges to a point $z^{\star}\in B(z;r)$ or the set of cluster
points of $(z_{n})_{n\in\mathbf{N}}$ forms a nontrivial continuum
in $B(z;r)$.\emph{ This proves the first part of claim }(i)\emph{
and claim }(ii)\emph{. }

\textbf{Step 4. }Now we prove the second part of claim (i). Under
the additional condition $\mathop{\underline{\mathrm{lim}}}_{n\to\infty}\sigma_{z_{n}z^{\star}}^{2}=0$,
we show that $z_{n}-\tilde{R}z_{n}\to0,z_{n}\to z^{\star}$ implies
$z^{\star}\in\mathop{{\bf fix}}\tilde{T}$, $\mathbf{prox}_{\gamma f}\left(z^{\star}\right)$
is an optimal solution of \eqref{eq:original_problem}, and $x_{n}\to\mathbf{prox}_{\gamma f}\left(z^{\star}\right)$.
For every $n\in\mathbf{N}$,
\begin{align}
\|z^{\star}-\tilde{R}z^{\star}\|^{2} & \overset{a)}{=}\|z_{n}-\tilde{R}z^{\star}\|^{2}-\|z_{n}-z^{\star}\|^{2}\nonumber \\
 & \quad-2\langle z_{n}-z^{\star}\mid z^{\star}-\tilde{R}z^{\star}\rangle\nonumber \\
 & =\|\big(z_{n}-\tilde{R}z_{n}\big)+\big(\tilde{R}z_{n}-\tilde{R}z^{\star}\big)\|^{2}\nonumber \\
 & \quad-\|z_{n}-z^{\star}\|^{2}-2\langle z_{n}-z^{\star}\mid z^{\star}-\tilde{R}z^{\star}\rangle\nonumber \\
 & =\|z_{n}-\tilde{R}z_{n}\|^{2}+\|\tilde{R}z_{n}-\tilde{R}z^{\star}\|^{2}\nonumber \\
 & \quad+2\langle z_{n}-\tilde{R}z_{n}\mid\tilde{R}z_{n}-\tilde{R}z^{\star}\rangle\nonumber \\
 & \quad-\|z_{n}-z^{\star}\|^{2}-2\langle z_{n}-z^{\star}\mid z^{\star}-\tilde{R}z^{\star}\rangle\nonumber \\
 & \overset{b)}{\leq}\|z_{n}-\tilde{R}z_{n}\|^{2}+\cancel{\|z_{n}-z^{\star}\|^{2}}+\sigma_{z_{n}z^{\star}}^{2}\nonumber \\
 & \quad+2\langle z_{n}-\tilde{R}z_{n}\mid\tilde{R}z_{n}-\tilde{R}z^{\star}\rangle-\cancel{\|z_{n}-z^{\star}\|^{2}}\nonumber \\
 & \quad-2\langle z_{n}-z^{\star}\mid z^{\star}-\tilde{R}z^{\star}\rangle\nonumber \\
 & =\|z_{n}-\tilde{R}z_{n}\|^{2}+2\langle z_{n}-\tilde{R}z_{n}\mid\tilde{R}z_{n}-\tilde{R}z^{\star}\rangle\nonumber \\
 & \quad-2\langle z_{n}-z^{\star}\mid z^{\star}-\tilde{R}z^{\star}\rangle+\sigma_{z_{n}z^{\star}}^{2},\label{eq:feedback}
\end{align}
where, $a)$ uses the identity 
\begin{align*}
\|z_{n}-\tilde{R}z^{\star}\|^{2} & =\|(z_{n}-z^{\star})+(z^{\star}-\tilde{R}z^{\star})\|^{2}\\
 & =\|z_{n}-z^{\star}\|^{2}+\|z^{\star}-\tilde{R}z^{\star}\|^{2}\\
 & \quad+2\left\langle z_{n}-z^{\star}\mid z^{\star}-\tilde{R}z^{\star}\right\rangle ,
\end{align*}
and $b)$ uses \eqref{eq:M_x_y}. We now compute the limit (or the
limit inferior) for each of the terms on the right-hand side of \eqref{eq:feedback}.
As $z_{n}-\tilde{R}z_{n}\to0$ and $z_{n}\to z^{\star}$, subtracting
them we have $\tilde{R}z_{n}\to z^{\star}$, hence $\tilde{R}z_{n}-\tilde{R}z^{\star}\to z^{\star}-\tilde{R}z^{\star}.$
Combining the last statement with $z_{n}-\tilde{R}z_{n}\to0$, we
have $\langle z_{n}-\tilde{R}z_{n}\mid\tilde{R}z_{n}-\tilde{R}z^{\star}\rangle\to0$.
Also, $z_{n}-z^{\star}\to0$ implies $\langle z_{n}-z^{\star}\mid z^{\star}-\tilde{R}z^{\star}\rangle\to0.$
Additionally, $\mathop{\underline{\mathrm{lim}}}_{n\to\infty}\sigma_{z_{n}z^{\star}}^{2}=0$.
So, using Lemma \ref{lem:-Let-df}, limit inferior of the right hand
side \eqref{eq:feedback} goes to zero. Hence, we conclude that $z^{\star}-\tilde{R}z^{\star}=0,$
\emph{i.e.,} $z^{\star}\in\mathop{{\bf fix}}\tilde{R}$. But, $\mathop{{\bf fix}}\tilde{R}=\mathop{{\bf fix}}\tilde{T}$
from Proposition \ref{prop:(Relationship-between-the}(ii). So, $z^{\star}\in\mathop{{\bf fix}}\tilde{T}$.
We now recall from Lemma \ref{lem:singtonness} that $\mathbf{prox}_{\gamma f}$
is continuous everywhere on $\mathbf{R}^{n}$. So, using the definition
of continuity, $z_{n}\to z^{\star}\in\mathop{{\bf fix}}\tilde{T}$
implies $x_{n+1}=\mathbf{prox}_{\gamma f}(z_{n})\to\mathbf{prox}_{\gamma f}(z^{\star})\in\mathbf{prox}_{\gamma f}(\mathop{{\bf fix}}\tilde{T})$.
But, $\mathbf{prox}_{\gamma f}(\mathop{{\bf fix}}\tilde{T})\subseteq\mathop{{\rm argmin}}(f+\delta_{\mathcal{C}})$
from Theorem \ref{prop:minimizer_charecterization}(iii). \emph{Thus
we have arrived at the second part of claim }(i)\emph{.}
\end{IEEEproof}

\subsubsection{Notes on Theorem \ref{(main-convergence-result)}\label{subsec:Notes-on-NC-DRS_conv}}

We make the following notes on Theorem \ref{(main-convergence-result)}.

$\bullet\,$\textbf{Nonemptiness of $\mathop{{\bf fix}}\tilde{T}.$}
Note that Theorem \ref{(main-convergence-result)} assumes that $\mathop{{\bf fix}}\tilde{T}$
is nonempty. This is a standard assumption in monotone operator theory
even in a convex setup \cite[\S 5.2]{Bauschke2017}.

$\bullet\,$\textbf{Relation to a convex setup. }In our convergence
analysis, the constraint set is nonempty and compact, but not necessarily
convex. However, our convergence analysis is also applicable to a
convex setup. Let $\mathcal{C}$ be convex. Then, both $(2\mathbf{\tilde{\Pi}}_{\mathcal{C}}-I_{n})$
and $(2\mathbf{prox}_{\gamma f}-I_{n})$ are nonexpansive operators,
hence, their composition $\tilde{S}=(2\mathbf{\tilde{\Pi}}_{\mathcal{C}}-I_{n})(2\mathbf{prox}_{\gamma f}-I_{n})$
is a nonexpansive operator. In such a convex setup, $\tilde{R}$ is
a nonexpansive operator, because $\tilde{S}=\tilde{R}$ from Proposition
\ref{prop:(Relationship-between-the} (where the relationship is established
irrespective of convexity). Then, recalling Remark \ref{(further-characterization-of},
expansiveness of $\tilde{R}$ is zero everywhere, \emph{i.e.,} $\sigma_{xy}=\varepsilon_{xy}=0$
at every $x,y$ in $\mathbf{R}^{n}$. As a result, the iteration scheme
\eqref{eq:convenient-form-1} corresponds to an averaged iteration
of a nonexpansive operator $\tilde{R}$, which guarantees convergence
of the sequence to a fixed point of $\tilde{R}$ for any initial point
\cite{ryu2016primer}. Also, the additional condition in the second
part of claim (i) are automatically satisfied. This guarantees the
convergence of \eqref{eq:DR_splitting_orig} to an optimal solution
for any initial point if we assume that $\mathcal{C}$ is convex. 

$\bullet\,$\textbf{Comments on the conditions.} Once we move from
a convex setup to a nonconvex setup, $\tilde{R}$ is not nonexpansive
anymore (recall Remark \ref{Cayley_is_not_nonexpansive}). Roughly
speaking, convergence in such a case requires that the total squared
expansiveness of $\tilde{R}$ stays bounded for the iterates with
respect to at least one fixed point of the nonconvex Douglas-Rachford
operator. More precisely, $\sum_{n=0}^{\infty}\sigma_{z_{n}z}^{2}$
needs to be bounded, where the sum represents the total deviation
of $\tilde{R}$ from being a nonexpansive operator over the sequence
$\left\{ \left(z_{n},z\right)\right\} _{n\in\mathbf{N}}$. If the
stated condition is satisfied, then $(z_{n})_{n\in\mathbf{N}}$ is
bounded in $B(z;r)$ and $\|z_{n+1}-z_{n}\|\to0$, but it does not
necessarily guarantee convergence to a point due to the lack of nonexpansiveness
of $\tilde{R}$, and this is why the cluster points of $(z_{n})_{n\in\mathbf{N}}$
may form a nontrivial continuum in $B(z;r)$.

Suppose now that $(z_{n})_{n\in\mathbf{N}}$ converges to a point
$z^{\star}$. Whether $z^{\star}$ is related to an optimal solution
of \eqref{eq:original_problem} would depend on $\mathop{\underline{\mathrm{lim}}}_{n\to\infty}\sigma_{z_{n}z^{\star}}^{2}$.
If it is zero, then $\mathbf{prox}_{\gamma f}\left(z^{\star}\right)$
is an optimal solution, and the iterate $x_{n}$ in \eqref{eq:DR_splitting_orig}
converges to this optimal solution. Roughly speaking, $\mathop{\underline{\mathrm{lim}}}_{n\to\infty}\sigma_{z_{n}z^{\star}}^{2}=0$
means that over $\left\{ \left(z_{n},z^{\star}\right)\right\} _{n\in\mathbf{N}}$,
$\tilde{R}$ acts as a nonexpansive operator in the lower limit. 

\section{Construction and convergence of \eqref{eq:ADMM_orig-1}}

In this section, we discuss how \eqref{eq:ADMM_orig-1} can be constructed
from the Douglas-Rachford splitting algorithm and comment on how the
construction influences the convergence properties of the former.
First, in \S\ref{subsec:Preliminaries_NC_ADMM} we present some preliminary
results to be used later. Then, in \S\ref{subsec:NC-ADMM-from-Douglas-Rachford}
we describe how \eqref{eq:ADMM_orig-1} is constructed from the Douglas-Rachford
splitting algorithm. Finally, in \S\ref{subsec:Convergence-properties-of-NC-ADMM}
we comment on convergence of \eqref{eq:ADMM_orig-1}, and we compare
it with \eqref{eq:DR_splitting_orig}.

\subsection{Preliminaries\label{subsec:Preliminaries_NC_ADMM} }

First, we describe the Douglas-Rachford splitting algorithm for minimizing
sum of two CPC functions; we will use it in the first step of constructing
\eqref{eq:ADMM_orig-1}. Then, we review the necessary background
on conjugate and biconjugate functions, and we present two lemmas
to be referenced in the second step of constructing \eqref{eq:ADMM_orig-1}.

\subsubsection{Douglas-Rachford splitting algorithm for minimizing sum of two CPC
functions}

Consider the convex optimization problem 
\begin{equation}
\begin{array}{ll}
\textup{minimize} & g\left(x\right)+h(x),\end{array}\label{eq:cvx-opt}
\end{equation}
where both $g:\mathbf{R}^{n}\to\overline{\mathbf{R}}$ and $h:\mathbf{R}^{n}\to\overline{\mathbf{R}}$
are CPC functions, and $x\in\mathbf{R}^{n}$ is the optimization variable.
The Douglas-Rachford splitting algorithm for this problem is
\begin{align}
x_{n+1} & =\mathbf{prox}_{\gamma h}\left(z_{n}\right)\nonumber \\
y_{n+1} & =\mathbf{prox}_{\gamma g}\left(2x_{n+1}-z_{n}\right)\tag{Convex-DRS}\label{eq:classic_DRS}\\
z_{n+1} & =z_{n}+y_{n+1}-x_{n+1},\nonumber 
\end{align}
where $n$ is the iteration counter, and $\gamma$ is a positive parameter.
In this convex setup, both $x_{n}$ and $y_{n}$ converge to an optimal
solution of \eqref{eq:cvx-opt} for any initial point \cite[Corollary 27.4]{Bauschke2017}.

\subsubsection{Conjugate and biconjugate of a function\label{subsec:Conjugate-and-biconjugate}}

Let $g:\mathbf{R}^{n}\to\left\{ -\infty\right\} \cup\overline{\mathbf{R}}$.
The \textbf{conjugate} of $g$, denoted by $g^{\star}$, is defined
as $g^{\star}\left(y\right)=\sup_{x\in\mathbf{R}^{n}}\left(\left\langle x\mid y\right\rangle -g\left(x\right)\right),$
which is closed and convex irrespective of the convexity of $g$ \cite[Proposition 13.11]{Bauschke2017}.
Also, the conjugate of a CPC function is CPC \cite[Theorem 4.3, Theorem 4.5]{beck2017first}.
Similarly, the \textbf{biconjugate} of $g$, denoted by $g^{\star\star}$,
is defined as $g^{\star\star}\left(y\right)=\sup_{x\in\mathbf{R}^{n}}\left(\left\langle x\mid y\right\rangle -g^{\star}\left(x\right)\right).$
Additionally, if the function is CPC, then its biconjugate is equal
to the function itself \cite[Lemma 4.8]{beck2017first}. Finally,
the relationship between the proximal operator of a CPC function $f$
with the proximal operator of its conjugate is given by \textbf{Moreau's
decomposition}: $\mathbf{prox}_{f}\left(x\right)+\mathbf{prox}_{f^{\star}}\left(x\right)=x$
for every $x\in\mathbf{R}^{n}$. Moreau's decomposition does not hold
for a nonconvex function.

Next, we present the following lemmas about conjugate functions in
the context of \eqref{eq:original_problem}. Here we use the notation
$g^{\vee},$ which denotes the reversal of a function $g$, and it
is defined as $g^{\vee}\left(x\right)=g\left(-x\right)$ for every
$x\in\mathbf{R}^{n}$. 
\begin{lem}[\textbf{proximal operator of $f^{\star\vee}$}]
\label{prox_gamma_f_star_vee} Let $f:\mathbf{R}^{n}\to\overline{\mathbf{R}}$
be the cost function in \eqref{eq:original_problem}. Then, for every
$\gamma>0$ and for every $x\in\mathbf{R}^{n}$, 
\[
\mathbf{prox}_{\gamma f^{\star\vee}}\left(x\right)=x+\gamma\,\mathbf{prox}_{\gamma^{-1}f}\left(-\gamma^{-1}x\right).
\]
\end{lem}
\begin{IEEEproof}
Recall that $f$ is CPC. For every $\gamma>0$ and for every $x\in\mathbf{R}^{n}$,
\begin{align*}
\mathbf{prox}_{\gamma f^{\star\vee}}\left(x\right) & \overset{a)}{=}\mathbf{prox}_{\gamma\left(f^{\vee}\right)^{\star}}\left(x\right)\\
 & \overset{b)}{=}x-\gamma\,\mathbf{prox}_{\gamma^{-1}f^{\vee}}\left(\gamma^{-1}x\right)\\
 & \overset{c)}{=}x+\gamma\,\mathbf{prox}_{\gamma^{-1}f}\left(-\gamma^{-1}x\right),
\end{align*}
where $a)$ follows from $f^{\star\vee}=f^{\vee\star}$ \cite[Proposition 13.20(v)]{Bauschke2017},
$b)$ follows from \cite[Proposition 23.29(viii)]{Bauschke2017} and
the fact that $f^{\vee}$ is CPC, and $c)$ directly follows from
\cite[Proposition 23.29(v)]{Bauschke2017}.
\end{IEEEproof}
In the following Lemma \textbf{convex hull} of a nonempty set $\mathcal{C}$,
which is the smallest convex set containing $\mathcal{C}$, is denoted
by $\mathop{{\bf conv}}\mathcal{C}$. Closure of $\mathop{{\bf conv}}\mathcal{C}$
is denoted by $\overline{\mathop{{\bf conv}}}\mathcal{C}$.

\begin{lem}[\textbf{conjugate and biconjugate of indicator function of $\mathcal{C}$}]
\emph{\label{biconjugate-of-indicator-of-C}} Let $\mathcal{C}$
be the constraint set in \eqref{eq:original_problem}. Then, 

\emph{(i)} $\delta_{\mathcal{C}}^{\star\star}=\delta_{\mathop{{\bf conv}}\mathcal{C}}$,
and 

\emph{(ii)} $\mathbf{prox}_{\gamma\delta_{\mathcal{C}}^{\star}}\left(x\right)=x-\gamma\mathbf{\Pi}_{\mathop{{\bf conv}}\mathcal{C}}\left(\gamma^{-1}x\right).$
\end{lem}
\begin{IEEEproof}
(i): From \cite[Example 4.2, Example 4.9]{beck2017first}, we have
$\delta_{\mathcal{C}}^{\star\star}=\delta_{\overline{\mathop{{\bf conv}}}\mathcal{C}}$.
The constraint set $\mathcal{C}$ is compact, hence its convex hull
$\mathop{{\bf conv}}\mathcal{C}$ is also compact, hence closed \cite[Corollary 2.30]{Rockafellar2009}.
So, $\overline{\mathop{{\bf conv}}}\mathcal{C}=\mathop{{\bf conv}}\mathcal{C}$,
and we conclude that $\delta_{\mathcal{C}}^{\star\star}=\delta_{\mathop{{\bf conv}}\mathcal{C}}$. 

(ii): As the constraint set $\mathcal{C}$ is nonempty and compact,
its indicator function $\delta_{\mathcal{C}}$ is closed \cite[Example 1.25]{Bauschke2017}
and proper \cite[page 7]{Rockafellar2009}. Hence, its conjugate $\delta_{\mathcal{C}}^{\star}$,
which is called the \textbf{support function} of the set $\mathcal{C}$,
is CPC (closed and convex due to \cite[Proposition 13.11]{Bauschke2017},
proper because $\mathcal{C}$ is bounded). As the conjugate of a CPC
function is CPC \cite[Theorem 4.3, Theorem 4.5]{beck2017first}, the
function $\delta_{\mathcal{C}}^{\star\star}$ is CPC. Using Moreau's
decomposition for every $x\in\mathbf{R}^{n}$,
\begin{align*}
\mathbf{prox}_{\gamma\delta_{\mathcal{C}}^{\star}}\left(x\right) & =x-\mathbf{prox}_{\left(\gamma\delta_{\mathcal{C}}^{\star}\right)^{\star}}\left(x\right)\\
 & \overset{a)}{=}x-\mathbf{prox}_{\gamma\delta_{\mathcal{C}}^{\star\star}\big(\gamma^{-1}(\cdot)\big)}\left(x\right)\\
 & \overset{b)}{=}x-\gamma\,\mathbf{prox}_{\gamma^{-1}\delta_{\mathcal{C}}^{\star\star}}\left(\gamma^{-1}x\right)\\
 & \overset{c)}{=}x-\gamma\mathbf{\Pi}_{\mathop{{\bf conv}}\mathcal{C}}\left(\gamma^{-1}x\right),
\end{align*}
where $a)$ follows from \cite[Proposition 13.20(i)]{Bauschke2017},
$b)$ follows from \cite[Proposition 23.29(iii)]{Bauschke2017}, and
$c)$ follows from combining $\delta_{\mathcal{C}}^{\star\star}=\delta_{\mathop{{\bf conv}}\mathcal{C}}$
in (i) and Corollary \ref{lem:Proximal-operator-of}.
\end{IEEEproof}

\subsection{Constructing \eqref{eq:ADMM_orig-1} from Douglas-Rachford splitting\label{subsec:NC-ADMM-from-Douglas-Rachford}}

This subsection is organized as follows. First, by applying \eqref{eq:classic_DRS}
to the convex dual of \eqref{eq:original_problem} we construct a
relaxed version of \eqref{eq:ADMM_orig-1}, where the projection is
onto $\mathop{{\bf conv}}\mathcal{C}$ rather than $\mathcal{C}$.
Then, we show that the relaxed version \eqref{eq:ADMM_orig-1} minimizes
$f$ over $\mathop{{\bf conv}}\mathcal{C}$. Next, we discuss construction
of \eqref{eq:ADMM_orig-1} from the relaxed variant by restricting
the latter's projection step onto $\mathcal{C}$. Finally, we comment
on the convergence properties of \eqref{eq:ADMM_orig-1} and relate
it to \eqref{eq:DR_splitting_orig}.

\subsubsection{Constructing dual of \eqref{eq:original_problem}}

Using indicator function, we write \eqref{eq:original_problem} as
\[
\begin{array}{ll}
\textup{minimize} & f\left(x\right)+\delta_{\mathcal{C}}\left(y\right)\\
\textup{subject to} & x-y=0,
\end{array}
\]
where $x,y\in\mathbf{R}^{n}$ are the optimization variables. Denote
the optimal value of the problem above by $p^{\star}$. The dual of
the reformulated problem, which is a convex optimization problem \cite[\S 5.1.6]{boyd2004convex},
is 
\[
\begin{array}{ll}
\textup{maximize} & -f^{\star\vee}\left(\nu\right)-\delta_{\mathcal{C}}^{\star}\left(\nu\right),\end{array}
\]
where $\nu\in\mathbf{R}^{n}$ is the optimization variable. Denote
the optimal value of the dual problem by $d^{\star}$. Due to weak
duality, we have, $d^{\star}\leq p^{\star}$, and, as the primal problem
is nonconvex, the duality gap $p^{\star}-d^{\star}$ is strict in
general. For convenience, we write the dual problem in minimization
form: 
\begin{equation}
\begin{array}{ll}
\textup{minimize} & f^{\star\vee}\left(\nu\right)+\delta_{\mathcal{C}}^{\star}\left(\nu\right),\end{array}\tag{Dual-OPT}\label{eq:dual-opt}
\end{equation}
with optimal value $-d^{\star}$ and same set of optimal solutions
as the dual problem. As $f$ is CPC, $f^{\star\vee}$ is also CPC
(follows from \S\ref{subsec:Conjugate-and-biconjugate} and \cite[Proposition 8.20]{Bauschke2017}).
Furthermore, from Lemma \ref{biconjugate-of-indicator-of-C}(ii),
$\delta_{\mathcal{C}}^{\star}$ is also CPC. Thus we can apply \eqref{eq:classic_DRS}
to \eqref{eq:dual-opt}. 

\subsubsection{Applying Douglas-Rachford splitting to \eqref{eq:dual-opt}}

By setting $g:=f^{\star\vee}$ and $h:=\delta_{\mathcal{C}}^{\star}$
in \eqref{eq:cvx-opt}, we have the following Douglas-Rachford splitting
algorithm for the dual problem:
\begin{align}
\zeta_{n+1} & =\mathbf{prox}_{\gamma\delta_{\mathcal{C}}^{\star}}\left(\psi_{n}\right)\nonumber \\
\xi_{n+1} & =\mathbf{prox}_{\gamma f^{\star\vee}}\left(2\zeta_{n+1}-\psi_{n}\right)\tag{Dual-DRS}\label{eq:dual-DRS}\\
\psi_{n+1} & =\psi_{n}+\xi_{n+1}-\zeta_{n+1}.\nonumber 
\end{align}

Using Lemma \ref{biconjugate-of-indicator-of-C} and Lemma \ref{prox_gamma_f_star_vee},
we simplify the first two iterates of \eqref{eq:dual-DRS} as
\begin{align*}
\zeta_{n+1} & =\psi_{n}-\gamma\,\mathbf{\Pi}_{\mathop{{\bf conv}}\mathcal{C}}\big(\gamma^{-1}\psi_{n}\big),\:\textrm{and}\\
\xi_{n+1} & =2\zeta_{n+1}-\psi_{n}+\gamma\,\mathbf{prox}_{\gamma^{-1}f}\big(-\gamma^{-1}(2\zeta_{n+1}-\psi_{n})\big).
\end{align*}

Using these simplified iterates and introducing intermediate iterates
$\tilde{y}_{n+1}=\mathbf{\Pi}_{\mathop{{\bf conv}}\mathcal{C}}\left(\gamma^{-1}\psi_{n}\right)$
and $\tilde{x}_{n+1}=\mathbf{prox}_{\gamma^{-1}f}\big(-\gamma^{-1}(2\zeta_{n+1}-\psi_{n})\big)$,
we can write \eqref{eq:dual-DRS} as
\begin{align*}
\tilde{y}_{n+1} & =\mathbf{\Pi}_{\mathop{{\bf conv}}\mathcal{C}}\big(\gamma^{-1}\psi_{n}\big)\\
\zeta_{n+1} & =\psi_{n}-\gamma\tilde{y}_{n+1}\\
\tilde{x}_{n+1} & =\mathbf{prox}_{\gamma^{-1}f}\big(-\gamma^{-1}(2\zeta_{n+1}-\psi_{n})\big)\\
 & =\mathbf{prox}_{\gamma^{-1}f}\big(-\gamma^{-1}(\psi_{n}-2\gamma\tilde{y}_{n+1})\big)\\
\xi_{n+1} & =2\zeta_{n+1}-\psi_{n}+\gamma\tilde{x}_{n+1}\\
 & =\psi_{n}-2\gamma\tilde{y}_{n+1}+\gamma\tilde{x}_{n+1}\\
\psi_{n+1} & =\psi_{n}+\xi_{n+1}-\zeta_{n+1}\\
 & =\psi_{n}-\gamma\tilde{y}_{n+1}+\gamma\tilde{x}_{n+1}.
\end{align*}

Note that the iterates $\zeta_{n}$ and $\xi_{n}$ do not have any
explicit dependence, hence they can be removed. Furthermore, introduce
a new iterate, $z_{n}=\frac{1}{\gamma}\psi_{n}-\tilde{x}_{n}$ . Substituting
$\psi_{n}:=\gamma\left(z_{n}+\tilde{x}_{n}\right)$ in the iteration
scheme above, we get
\begin{align}
\tilde{y}_{n+1} & =\mathbf{\Pi}_{\mathop{{\bf conv}}\mathcal{C}}\left(z_{n}+\tilde{x}_{n}\right)\nonumber \\
\tilde{x}_{n+1} & =\mathbf{prox}_{\gamma^{-1}f}\left(-\left(z_{n}+\tilde{x}_{n}-2\tilde{y}_{n+1}\right)\right)\nonumber \\
 & \overset{a)}{=}\mathbf{prox}_{\gamma^{-1}f}\left(-\left(z_{n+1}-\tilde{y}_{n+1}\right)\right)\nonumber \\
z_{n+1} & =z_{n}+\tilde{x}_{n}-\tilde{y}_{n+1},\label{eq:idw}
\end{align}
where $a)$ follows from \eqref{eq:idw}.

Finally, we swap the order of of $\tilde{x}_{n+1}$ and $z_{n+1}$
to obtain the correct dependency:
\begin{align*}
\tilde{y}_{n+1} & =\mathbf{\Pi}_{\mathop{{\bf conv}}\mathcal{C}}\left(z_{n}+\tilde{x}_{n}\right)\\
z_{n+1} & =z_{n}+\tilde{x}_{n}-\tilde{y}_{n+1},\\
\tilde{x}_{n+1} & =\mathbf{prox}_{\gamma^{-1}f}\left(\tilde{y}_{n+1}-z_{n+1}\right).
\end{align*}
We now substitute $\tilde{x}_{n}:=x_{n+1}$, $\tilde{y}_{n}:=y_{n}$,
and $\frac{1}{\gamma}:=\tilde{\gamma}$ in the iterates above to obtain
a relaxed version \eqref{eq:ADMM_orig-1}:
\begin{align}
x_{n+1} & =\mathbf{prox}_{\tilde{\gamma}f}\left(y_{n}-z_{n}\right)\nonumber \\
y_{n+1} & =\mathbf{\Pi}_{\mathop{{\bf conv}}\mathcal{C}}\left(z_{n}+x_{n+1}\right)\tag{Relaxed-NC-ADMM}\label{eq:ADMM-to-dual}\\
z_{n+1} & =z_{n}-y_{n+1}+x_{n+1},\nonumber 
\end{align}
which is similar to \eqref{eq:ADMM_orig-1}, except the projection
is onto $\mathop{{\bf conv}}\mathcal{C}$ rather than onto $\mathcal{C}$. 

\subsubsection{Constructing \eqref{eq:ADMM_orig-1} from \eqref{eq:ADMM-to-dual}}

Now we discuss how we can arrive at \eqref{eq:ADMM_orig-1} from \eqref{eq:ADMM-to-dual}.
The first step requires the observation that \eqref{eq:ADMM-to-dual}
finds a minimizer of $f$ over $\mathop{{\bf conv}}\mathcal{C}$.
To see that, construct the dual of \eqref{eq:dual-opt}, which is
\begin{equation}
\begin{array}{ll}
\textup{maximize} & -\left(f^{\star\vee}\right)^{\star\vee}\left(x\right)-\delta_{\mathcal{C}}^{\star\star}\left(x\right),\end{array}\tag{Double-Dual}\label{eq:double-dual}
\end{equation}
where $x\in\mathbf{R}^{n}$ is the optimization variable. As both
\eqref{eq:dual-opt} and \eqref{eq:double-dual} are convex optimization
problems, strong duality usually holds (under constraint qualifications),
where both problems have the same optimal value $-d^{\star}$. Now,
$\left(f^{\star\vee}\right)^{\star\vee}\overset{a)}{=}\left(f^{\vee\star}\right)^{\star\vee}\overset{}{=}\left(\left(f^{\vee}\right)^{\star\star}\right)^{\vee}\overset{b)}{=}f^{\vee\vee}\overset{c)}{=}f,$where
$a)$ follows from $f^{\star\vee}=f^{\vee\star}$ for CPC function
$f$ \cite[Proposition 13.20(v)]{Bauschke2017}, $b)$ follows from
the fact that the biconjugate of a CPC function is equal to the function
itself \cite[Lemma 4.8]{beck2017first}, and $c)$ follows from the
fact that applying reversal operation twice on a function returns
the original function. Furthermore, $\delta_{\mathcal{C}}^{\star\star}=\delta_{\mathop{{\bf conv}}\mathcal{C}}$
from Lemma \ref{biconjugate-of-indicator-of-C}(i). Hence, the dual
of \eqref{eq:double-dual}, written as a minimization problem, is
\[
\begin{array}{ll}
\textup{minimize} & f\left(x\right)\\
\textup{subject to} & x\in\mathop{{\bf conv}}\mathcal{C},
\end{array}
\]
where $x\in\mathbf{R}^{n}$ is the optimization variable with optimal
value $d^{\star}$. So, under strong duality between \eqref{eq:dual-opt}
and \eqref{eq:double-dual}, \eqref{eq:ADMM-to-dual} finds a minimizer
of $f$ over the set $\mathop{{\bf conv}}\mathcal{C}$, which appears
in the projection step of \eqref{eq:ADMM-to-dual}. So, to solve the
original problem \eqref{eq:original_problem} where we seek a minimizer
of $f$ over $\mathcal{C}$, an intuitive modification to \eqref{eq:ADMM-to-dual}
is replacing $\mathop{{\bf conv}}\mathcal{C}$ with $\mathcal{C}$
(hence $\mathbf{\Pi}_{\mathop{{\bf conv}}\mathcal{C}}$ with $\tilde{\mathbf{\Pi}}_{\mathcal{C}}$),
which results in \eqref{eq:ADMM_orig-1}. Roughly speaking, \eqref{eq:ADMM_orig-1}
is constructed by first relaxing the constraint set of the original
problem to its convex hull, then applying the Douglas-Rachford splitting
algorithm for the relaxed problem and finally restricting the resultant
algorithm on the original constraint set. 

The construction procedure also provides an alternative explanation
behind why, when compared with exact solvers, \eqref{eq:ADMM_orig-1}
often achieves lower objective values in many numerical experiments
performed in \cite{Diamond2018,Takapoui2017a,Takapoui2017,takapoui2016simple}.
In these works, these lower objective values are attributed to the
superior performance of \eqref{eq:ADMM_orig-1} in solving nonconvex
problems based on empirical evidence. An alternative explanation could
be that the heuristic is solving a modified dual problem, which, in
the absence of strong duality, is guaranteed to yield an objective
value that is smaller than or equal to that of the original problem.

\subsection{Convergence of \eqref{eq:ADMM_orig-1}\label{subsec:Convergence-properties-of-NC-ADMM}}

Now we comment on convergence properties of \eqref{eq:ADMM_orig-1}
in comparison with \eqref{eq:DR_splitting_orig}.

$\bullet\,$\textbf{Convergence to an optimal solution.} For \eqref{eq:DR_splitting_orig},
the fixed point set of $R_{\gamma\partial\delta_{\mathcal{C}}}R_{\gamma\partial f}$
acts as a bridge between global minimizers of \eqref{eq:original_problem}
and  the fixed point set of the nonconvex Douglas-Rachford operator
(Theorem \ref{prop:minimizer_charecterization}). Though \eqref{eq:classic_DRS}
is equivalent to \eqref{eq:ADMM-to-dual} under strong duality, no
such equivalence seems to exist between \eqref{eq:DR_splitting_orig}
and \eqref{eq:ADMM_orig-1}, because there is a strict duality gap
between \eqref{eq:original_problem} and \eqref{eq:dual-opt}, and
$\mathbf{\Pi}_{\mathop{{\bf conv}}\mathcal{C}}\neq\mathbf{\tilde{\Pi}}_{\mathcal{C}}$.
This lack of equivalence prevents connecting the fixed point set of
the underlying \eqref{eq:ADMM_orig-1} operator to global minimizers
of \eqref{eq:original_problem} through the fixed point set of $R_{\gamma\partial\delta_{\mathcal{C}}}R_{\gamma\partial f}$.

$\bullet\,$\textbf{Convergence to a point. }Furthermore, the lack
of equivalence between \eqref{eq:DR_splitting_orig} and \eqref{eq:ADMM_orig-1}
makes it harder to comment analogously on convergence of \eqref{eq:ADMM_orig-1}
to a general point (not necessarily an optimal solution) as well.
As shown in the proof of Theorem \ref{(main-convergence-result)},
establishing convergence to a general point for \eqref{eq:DR_splitting_orig}
depends on the interrelationship between the nonconvex Douglas-Rachford
operator and the nonconvex Peaceman-Rachford operator. Unfortunately,
such relationship may break down for \eqref{eq:ADMM_orig-1}, because
constructing such a relationship would require Moreau's decomposition
to hold for nonconvex functions.

\section{Future work}

Future research directions include conducting numerical experiments
to compare the performance of \eqref{eq:DR_splitting_orig} with \eqref{eq:ADMM_orig-1}.

\bibliographystyle{IEEEtran}
\bibliography{lattice_based_ADMM_manuscript}

\end{document}